\documentclass[12pt]{elsarticle}     

\makeatletter
\def\ps@pprintTitle{%
 \let\@oddhead\@empty
 \let\@evenhead\@empty
 \def\@oddfoot{}%
 \let\@evenfoot\@oddfoot}
\makeatother

\usepackage{amsmath,amscd,wrapfig,amsfonts,amssymb,mathtools,mathrsfs}
\usepackage{booktabs,multirow,lscape,cool,bm}
\usepackage{url,parskip,caption,subcaption}

\usepackage[top=1.10in, bottom=1.10in, left=1.10in, right=1.10in]{geometry}

\makeatletter
\g@addto@macro\normalsize{%
  \setlength\abovedisplayskip{.6em}
  \setlength\belowdisplayskip{.6em}
  \setlength\abovedisplayshortskip{.6em}
  \setlength\belowdisplayshortskip{.6em}
}

\begin{document}








\begin{frontmatter}

\begin{keyword}
fast algorithms\sep
ordinary differential equations
\end{keyword}

\title
{
A frequency-independent solver for  systems of first order linear ordinary differential equations
}

\begin{abstract}
When a system of first order linear ordinary differential equations has eigenvalues of large magnitude,
its solutions exhibit complicated behaviour, such as high-frequency oscillations,
rapid growth or rapid decay.  The cost of representing such solutions using standard
techniques typically grows with the magnitudes of the eigenvalues.
As a consequence, the running times of  standard solvers for ordinary differential equations
also grow with the size of these eigenvalues.
The solutions of  scalar equations with slowly-varying 
coefficients, however,  can be efficiently represented via slowly-varying phase functions, 
regardless of the magnitudes of the eigenvalues of the corresponding coefficient matrix.
Here, we couple an existing solver for scalar equations which exploits this observation 
with a well-known technique for transforming a system of linear ordinary differential equations into
scalar form.   The result is a method for solving a large class of systems of linear ordinary
differential equations in  time independent of the magnitudes of the eigenvalues of their coefficient
matrices.    We discuss the results of numerical experiments demonstrating the properties of our algorithm.

\end{abstract}


\author[1]{Tony Hu}
\author[2]{James Bremer\corref{cor1}}
\ead{bremer@math.toronto.edu}
\cortext[cor1]{Corresponding author}
\address[1]{University of Toronto}
\address[2]{Department of Mathematics,  University of Toronto}

\end{frontmatter}


\begin{section}{Introduction}

The complexity of the solutions of a system of first order linear ordinary differential equations
\begin{equation}
\mathbf{y}'(t) = A(t) \mathbf{y}(t)
\label{introduction:system}
\end{equation}
increases with the magnitudes of the eigenvalues $\lambda_1(t),\ldots,\lambda_n(t)$ of the coefficient 
matrix $A(t)$.  Indeed, the cost to represent these solutions over an interval $[a,b]$ using standard methods, such as polynomial
and trigonometric expansions, grows roughly linearly with the quantity
\begin{equation}
\Omega = \max_{\substack{i=1,2,\ldots,n}} \int_a^b \left| \lambda_i(t)\right|\, dt,
\label{introduction:omega}
\end{equation}
%
which we refer to  the ``frequency'' of the system (\ref{introduction:system}).  We choose
this  terminology  because, in most cases of interest, it is the imaginary parts of 
$\lambda_1(t),\ldots,\lambda_n(t)$ that are of large magnitude.  In fact,
when one or more of the eigenvalues has a real part of large magnitude,
initial and boundary value problems for (\ref{introduction:system}) are highly ill-conditioned
and solving them numerically requires additional information about the desired solution.

However, a large class of scalar ordinary differential equations of the form
\begin{equation}
y^{(n)}(t) + q_{n-1}(t) y^{(n-1)}(t) + \cdots + q_1(t) y'(t) + q_0(t) y(t) = 0
\label{introduction:scalarode}
\end{equation}
admit phase functions whose cost to represent via standard methods, such as polynomial or trigonometric
expansions, depends on the  complexity of the coefficients $q_0,\ldots,q_{n-1}$
but not the magnitudes of the eigenvalues  $\lambda_1(t), \ldots, \lambda_n(t)$ 
of the corresponding coefficient matrix
\begin{equation}
\left(\begin{array}{cccccccc}
0 & 1 & 0 & \cdots & 0 & 0 \\
0 & 0 & 1 & \cdots & 0 & 0 \\
\vdots &  \vdots& \vdots & \ddots & \vdots& \vdots\\
0 & 0 & 0 & \cdots & 1 & 0 \\
0 & 0 & 0 & \cdots & 0 & 1 \\
-q_0(t) & -q_1(t) & -q_2(t) & \cdots&  -q_{n-2}(t) & -q_{n-1}(t)\\
\end{array}
\right).
\label{introduction:scalarcoef}
\end{equation}
Indeed, if the $q_0,\ldots,q_{n-1}$ are slowly-varying on an interval $I$ and 
$\lambda_1(t), \ldots, \lambda_n(t)$ are
 distinct for all  $t$ in $I$, then it is possible to find slowly-varying functions 
$\psi_1,\ldots,\psi_n\colon I \to \mathbb{C}$
such that 
\begin{equation}
y_j(t) = \exp\left(\psi_j(t)\right), \ \ \ \ \ j=1,\ldots,n,
\label{introduction:phaserep}
\end{equation}
form a basis in the space of solutions of (\ref{introduction:scalarode}) given on the interval $I$.
That slowly-varying phase functions exist under these conditions, at least in an asymptotic sense,
has long been known and this observation is the basis of the WKB method and other related techniques 
(see, for instance, \cite{Miller}, \cite{Wasov} and \cite{SpiglerPhase1,SpiglerPhase2,SpiglerZeros}).  
A careful proof of the existence of slowly-varying phase functions for second order differential
equations, which can be extended to higher order scalar equations, is given in \cite{BremerRokhlin}.  
We say that the  system (\ref{introduction:system}) is nondegenerate on the interval $I$
if the condition mentioned above holds.  That is, provided all of the eigenvalues $\lambda_1(t), \ldots, \lambda_n(t)$  of its 
coefficient matrix $A(t)$ are distinct for each $t \in I$.  Moreover, we say that $t_0$ is a turning point for 
(\ref{introduction:system}) if the eigenvalues  are distinct in a deleted neighbourhood of $t_0$ but coalesce at $t_0$.

It is only relatively recently that numerical algorithms which exploit the existence
of slowly-varying phase functions have appeared.  The algorithm of \cite{BremerPhase} uses a continuation  scheme of sorts
to  construct slowly-varying phase functions $\psi_1$ and $\psi_2$ such that
\begin{equation}
\exp(\psi_1(t)) \ \ \mbox{and} \ \ \ \exp\left(\psi_2(t)\right)
\end{equation}
form a basis in the space of solutions of a second order differential equation of the form
\begin{equation}
y''(t) + q(t) y(t) = 0, \ \ \ \ \ a \leq t \leq b,
\label{introduction:two}
\end{equation}
where $q$ is slowly-varying and does not change sign on the interval $(a,b)$.   Under mild conditions on $q$, 
the running time of this algorithm is independent of the size of the eigenvalues of the system corresponding
to (\ref{introduction:two}), which are, of course, given by 
\begin{equation}
\lambda_1(t) = \sqrt{-q(t)}\ \ \ \mbox{and} \ \ \ \lambda_2(t) = - \sqrt{-q(t)}.
\label{introduction:twoeigs}
\end{equation}
The condition that $q$ does not change sign implies that  (\ref{introduction:two}) is nondegenerate on $(a,b)$.
This condition is imposed because slowly-varying phase function need not extend across
turning points.  The  method can, however, be extended to the case in which  (\ref{introduction:two})
is nondegenerate on  $[a,b]$ except  at a finite number of turning points.
Indeed,  if  $a=\xi_1 < \xi_2 < \ldots < \xi_k=b$ is a partition
of $[a,b]$ such that $\xi_2,\ldots,\xi_{k-1}$ are the roots of $q$ in the open
interval $(a,b)$, then simply applying the method of \cite{BremerPhase}
to each of the subintervals $\left[\xi_j,\xi_{j+1}\right]$, $j=1,\ldots,k-1$,  results in a collection
of $2(k-1)$ slowly-varying phase functions which efficiently represent the solutions
of (\ref{introduction:two}).  A detailed discussion of this approach, including
an account of many numerical experiments performed to demonstrate its efficacy, 
can be found in \cite{BremerPhase2}.

The continuation method of \cite{BremerPhase} readily extends to the case of
higher order scalar equations with slowly-varying coefficients.
However, the authors have found that a strongly related class of algorithms introduced in \cite{Scalar1} performs 
slightly better in practice.  When applied to an equation of the form (\ref{introduction:scalarode})
with slowly-varying coefficients, the  methods of \cite{Scalar1} produce a collection
$\psi_1,\ldots,\psi_n$ of slowly-varying phase functions such that
(\ref{introduction:phaserep})
is a basis in the space of  solutions of (\ref{introduction:scalarode}). 
The running times of these algorithms are  largely independent of the magnitude
of the eigenvalues of the corresponding coefficient matrix (\ref{introduction:scalarcoef}).

Here, we couple one of the  algorithms of \cite{Scalar1} 
with a standard technique for reducing systems of linear ordinary 
differential equations to scalar form.  The result is a solver for 
initial and boundary value problems for a large class of linear systems of differential
equations which runs in time independent of frequency.
As in \cite{BremerPhase}, we focus here on the case in which
the system is nondegenerate on the interior of the domain $[a,b]$ of the differential equation.  However, our algorithm
can  easily be extended to systems which are nondegenerate on an interval $[a,b]$ except at a finite
number of turning points by applying it repeatedly, on a collection of subintervals of $[a,b]$, as in \cite{BremerPhase2}.

There is, though, one significant difficulty with our method.  The transformation matrices we form
to convert the system (\ref{introduction:system}) to the scalar form  (\ref{introduction:scalarcoef})
can be ill-conditioned, and  the condition numbers of these transformation matrices
tend to increase with the dimension of the system under consideration.
As a consequence, the obtainable accuracy of our method also generally decreases as the dimension of the system
increases.
For the systems of two and three equations we considered (which include those discussed in the experiments of this paper as well as many  others we
considered while developing and testing our algorithm), we found this effect to be extremely mild.  Indeed, for these systems, 
we found that the solutions could easily be calculated to  around 12 digits of relative accuracy at low frequencies.  
The errors increased further as the frequency of the problems increased; however,  this is to be expected from any numerical
method as  the condition number of initial and boundary value problems for (\ref{introduction:system}) also grow with frequency.
For systems of four equations, around 10 digits of relative accuracy
could be readily  obtained at low frequencies.   Systems of dimension five and higher can be addressed by our algorithm,
but the obtainable accuracy continues to deteriorate.  

Despite the difficulties in constructing well-conditioned transformations which convert a system to scalar form, the algorithm of 
this paper applies to a large class of systems of linear ordinary differential equations of modest orders 
and our method appears to be the first high-accuracy solver for such problems which runs in time independent of frequency.
We view this work as a step toward developing robust frequency-independent solvers for 
systems of ordinary differential equations.    We discuss  several obvious directions for further development
 in the conclusion of this paper.

Because  the cost of standard solvers becomes prohibitive when $\Omega$ is large,
many specialized techniques have been developed for solving systems of linear
ordinary differential equations in this regime.  At the present time, the most widely-used 
such methods are based on Magnus expansions.
%
%
Introduced in \cite{Magnus},    Magnus expansions are certain series of the form
\begin{equation}
\sum_{k=1}^\infty \Omega_k(t)
\label{introduction:Magnus}
\end{equation}
such that  $\exp\left(\sum_{k=1}^\infty \Omega_k(t)\right)$
locally represents a fundamental matrix for a system of  differential equations
\begin{equation}
\mathbf{y}'(t) = A(t) \mathbf{y}(t).
\end{equation}
The first few terms for the series around $t=0$ are given by 
\begin{equation}
\begin{aligned}
\Omega_1(t) &=  \int_{0}^t A(s)\ ds,\\
\Omega_2(t) &= \frac{1}{2} \int_{0}^t \int_{0}^{t_1} \left[A(t_1),A(t_2)\right]\, dt_2dt_1 \ \ \ \mbox{and} \\
\Omega_3(t) &= \frac{1}{6} \int_{0}^t \int_{0}^{t_1} \int_{0}^{t_2} 
\left[A(t_1), \left[A(t_2),A(t_3)\right]\right] +  \left[A(t_3), \left[A(t_2),A(t_1)\right]\right]\,  dt_3dt_2dt_1.
\end{aligned}
\label{introduction:magnusterms}
\end{equation}
The straightforward evaluation of the $\Omega_j$ is nightmarishly expensive; however, a clever
technique which renders the calculations manageable is introduced in \cite{Iserles101}
and it paved the way for the development of a class of  numerical solvers 
which represent a fundamental matrix for 
(\ref{introduction:system}) over an interval $I$ via a collection of truncated
Magnus expansions. While the entries of the $\Omega_j$ are slowly-varying whenever the entries of $A(t)$ are slowly-varying, the 
radius of convergence of the series  in (\ref{introduction:Magnus})
depends on the magnitude of the coefficient matrix $A(t)$, which is, in turn,
related to the magnitudes of the eigenvalues of $A(t)$.
Of course, this means that the number of Magnus expansions
which are needed, and hence the cost of the method,
depends on the  magnitudes of the eigenvalues of $A(t)$.
See, for instance, \cite{Iserles102}, which gives for estimates
of the growth in the running time of Magnus expansion methods in 
the case of an equation of the form (\ref{introduction:two}) as 
a function of the magnitude of the coefficient $q$.
The algorithm of this paper, by contrast,
represents a  fundamental matrix for (\ref{introduction:system}) 
in the form 
\begin{equation}
\Phi(t)^{-1}
\left(
\begin{array}{cccccc}
u_1(t)          & u_2(t)         & \cdots            &   u_n(t)           & \\
u_1'(t)         & u_2'(t)        & \cdots            &   u_n'(t)            & \\
\vdots          & \vdots         & \ddots      &    \vdots           & \\
u_1^{(n-1)}(t)   & u_2^{(n-1)}(t)   & \cdots            &   u_n^{(n-1)}(t)             & \\
\end{array}
\right),
\label{introduction:expdiag}
\end{equation}
where the $u_j$ are given by 
\begin{equation}
u_j(t) = \exp(\psi_j(t)),\ \ \ \ \ j=1,\ldots,n,
\end{equation}
and the entries of $\Phi(t)^{-1}$ and  each of the $\psi_j$ are slowly-varying functions which  can be represented at a cost independent of the 
magnitudes of the eigenvalues of $A(t)$.


The remainder of this article is structured as follows.  In Section~\ref{section:reduction}, we discuss
the reduction of systems of linear ordinary differential equations to scalar form.
Section~\ref{section:scalar} discusses the approach of \cite{Scalar1} to the construction
of slowly-varying phase functions for scalar equations.
In Section~\ref{section:algorithm}, we detail the algorithm of this paper.
Section~\ref{section:experiments} discusses the results of numerical experiments conducted to 
demonstrate the properties of our algorithm.  We close with a few brief remarks
in Section~\ref{section:conclusions}.

\end{section}

\begin{section}{Reduction of a system of ordinary differential equations to a scalar equation}
\label{section:reduction}

It is well known that essentially any system of $n$ linear ordinary differential equations 
in $n$ unknowns can be  transformed into an $n^{th}$ order scalar  equation.  This is a standard result in 
differential Galois theory (see, for instance, \cite{DiffAlg} or \cite{PutSinger})  which goes back at least to 
\cite{Loewy}.  In this section, we discuss the well-known mechanism our algorithm uses for constructing such a transformation.

We first observe that if $I$ is an open subinterval of $\mathbb{R}$ and $\Phi \colon I \to \mathbb{C}$
is invertible for all $t \in I$, then letting
\begin{equation}
\mathbf{z}(t) = \Phi(t) \mathbf{y}(t)
\end{equation}
transforms  the system 
\begin{equation}
\mathbf{y}'(t) = A(t) \mathbf{y}(t) 
\end{equation}
into  
\begin{equation}
\mathbf{z}'(t) = B(t) \mathbf{z}(t),
\end{equation}
where
\begin{equation}
B(t) = \Phi'(t) \Phi(t)^{-1} + \Phi(t) A(t) \Phi(t)^{-1}.
\label{reduction:bt}
\end{equation}
Our goal is to construct $\Phi(t)$ in such a way that (\ref{reduction:bt}) is of the scalar form
\begin{equation}
\left(\begin{array}{cccccc}
0       & 1      & 0 & \cdots & 0 \\
0       & 0      & 1 & \cdots & 0 \\
\vdots  & \vdots & \vdots & \ddots & \vdots \\
0       & 0      & 0 & \cdots & 1 \\
-q_0(t) & -q_1(t) & -q_{2}(t) & \cdots & -q_{n-1}(t) \\
\end{array}\right).
\label{reduction:scalar}
\end{equation}
Transformations of this type correspond to cyclic vector for the transpose of $A(t)$,
which are smooth maps $\mathbf{v} \colon I \to \mathbb{C}^n$  such that the matrix
\begin{equation}
\Phi(t) = 
\left(\begin{array}{ccccc}
\mathbf{v}(t) \\
D\left[\mathbf{v}\right](t) \\
  D^2\left[\mathbf{v}\right](t)  \\
 \vdots \\
  D^{n-1}\left[\mathbf{v}\right](t)
\end{array}
\right),
\label{reduction:transform}
\end{equation}
where $D$ is the operator defined via
\begin{equation}
D\left[\mathbf{v}\right](t) = \mathbf{v}'(t) + (A(t))^t \mathbf{v}(t),
\end{equation}
is invertible for all $t \in I$  It is easy to see that when  $\Phi(t)$ is of the form (\ref{reduction:transform}),
the matrix
\begin{equation}
\Phi'(t) + \Phi(t) A(t)
\end{equation}
is given by 
\begin{equation}
\left(
\begin{array}{c}
D\left[\mathbf{v}\right](t) \\
D^2\left[\mathbf{v}\right](t) \\
\vdots\\
D^{n}\left[\mathbf{v}\right](t) \\
\end{array}
\right),
\end{equation}
so that  $B(t)$ 
is of  the scalar form (\ref{reduction:scalar}).

Because the set of all cyclic vectors is dense in a certain sense,
it is easy to find them.   Indeed, one of the standard mechanism used by computer algebra
systems to convert systems of differential equations to scalar form entails simply 
choosing random polynomial vector fields $\mathbf{v}(t)$ until a cyclic vector is found.
We refer the reader to \cite{DiffAlg} and its references for details.

The approach we use to find cyclic vectors is similar.  Our algorithm takes as input a 
constant vector $\mathbf{v}$  which generates the transformation $\Phi(t)$ via 
(\ref{reduction:transform}).   The user is expected to  make arbitrary choices of this vector
until  one produces a suitable transformation matrix.  It is important to note that it is not sufficient 
for $\Phi(t)$ to be merely invertible, as it is in symbolic computations.
The condition number of transformation matrix affects the obtainable accuracy
of the algorithm, and so various choices of $\mathbf{v}$ must be made until
a relatively well-conditioned transformation is found.   As mentioned
in the introduction, it appears to become increasingly difficulty to find a suitable
transformation matrix as the dimension of the system under consideration grows.

Many other mechanisms for converting systems of differential equations to scalar
form which have been previously proposed.    For instance, in \cite{Barkatou}, an 
procedure analagous to Gaussian elimination for transforming a system in to 
scalar form is introduced.  It applies a sequence of ``elementary
operations'' (which are somewhat more complicated than those used
in Gaussian elimination) to incrementally convert the coefficient
matrix  into the desired form.  A partially pivoted version of this algorithm aimed at
producing a numerically well-conditioned transformation matrix could almost
certainly be developed.  Moreover, it is clear that various well-known optimization
algorithms could be adapted to the problem of finding well-conditioned transformation matrices.

In the algorithm we present here, we use the simplest possible approach --- repeatedly make
arbitrary choices until a reasaonble outcome is obtained --- and we leave investigations into
improved methods for future work.

\end{section}

\begin{section}{Phase functions for scalar equations}
\label{section:scalar}

In this section, we discuss the construction  of slowly-varying
phase functions for scalar equations of the form (\ref{introduction:scalarode})
with slowly-varying coefficients.
One obvious method  entails solving the nonlinear $(n-1)^{st}$ order scalar  equation
satisfied by the derivatives $r_1,\ldots,r_n$ of the phase functions $\psi_1,\ldots,\psi_n$
which is obtained by  inserting the representation
\begin{equation}
y(t) = \exp\left(\int r(t)\, dt\right)
\label{reduction:log}
\end{equation}
into (\ref{introduction:scalarode}).   By a slight abuse
of terminology, we call this nonlinear equation the  $(n-1)^{st}$ order Riccati
equation, or the Riccati equation corresponding to (\ref{introduction:scalarode}).
The general form of this equation is quite complicated, but it relatively
simple at low orders.  When $n=2$, the equation is
\begin{equation}
r'(t) + (r(t))^2 + q_1(t) r(t) + q_0(t) = 0;
\label{scalar:riccati1}
\end{equation}
when $n=3$, it is
\begin{equation}
r''(t) + 3 r'(t) r(t) + (r(t))^3 + q_2(t) r'(t) + q_2(t) (r(t))^2 + q_1(t) r(t) + q_0(t) = 0;
\label{scalar:riccati2}
\end{equation}
and, for $n=4$, the equation is
\begin{equation}
\begin{aligned}
r'''(t) &+ 4 r''(t) r(t) + 3 (r'(t))^2 + 6  r'(t) (r(t))^2 + (r(t))^4 +q_3(t) (r(t))^3 + q_3(t) r''(t) \\
&+3 q_3(t) r'(t) r(t) +q_2(t) (r(t))^2 +q_2(t) r'(t) + q_1(t) r(t) +q_0(t) = 0.
\end{aligned}
\label{scalar:riccati3}
\end{equation}

Because most solutions of  Riccati equation corresponding to (\ref{introduction:scalarode}) 
are rapidly-varying when the eigenvalues $\lambda_1(t),\ldots,\lambda_n(t)$ of (\ref{introduction:scalarcoef})
are of large magnitude,  some mechanism is needed to  select the slowly-varying solutions.  
The algorithm of \cite{BremerPhase}, which applies to second order linear ordinary differential equations
of the form
\begin{equation}
y''(t) + q(t) y(t) = 0,\ \ \ \ a < t < b,
\label{scalar:two}
\end{equation}
with slowly-varying coefficients, uses a continuation scheme of sorts to do so.
More explicitly, it introduces a smoothly deformed version of the coefficient $\tilde{q}$ 
which equal to an appropriately chosen constant $\lambda^2$ on a small interval 
$[a_0,b_0]$ in $[a,b]$ and agrees with  $q$ outside of a neighbourhood of $U$ of $[a_0,b_0]$.  The values of two 
slowly-varying solutions  $\widetilde{r}_1$ and $\widetilde{r_2}$ of the Riccati equation 
corresponding to the deformed version 
\begin{equation}
y''(t) + \widetilde{q}(t) y(t) = 0,\ \ \ \ a < t < b,
\label{scalar:twodeformed}
\end{equation}
of (\ref{scalar:two}) are known at  $c$; indeed, we can take $\widetilde{r_1}(c) = i\lambda$  and $\widetilde{r_2}(c) = -i\lambda$.
Solving the Riccati equation corresponding to (\ref{scalar:twodeformed}) using these as initial
values allows us to compute the values of the derivatives of the derivatives of two slowly-varying
phase functions $r_1$ and $r_2$ for the original (\ref{scalar:two}) outside of the neighborhood $U$ of $[a_0,b_0]$.  We can then solve
the Riccati equation for the original equation to find the values of $r_1$ and $r_2$ inside $U$.
 This technique could be easily generalized to the case of higher order scalar equations, but the authors
have found the approach of  \cite{Scalar1}, which is inspired by the classical Levin scheme for
numerical evaluation of oscillatory integrals, to be somewhat more effective.

Introduced in  \cite{Levin}, the Levin method is based on the observation that inhomogeneous equations of the form
\begin{equation}
y'(t) + p_0(t) y(t)  = f(t)
\end{equation}
admit solutions whose complexity depends on that of $p_0$ and $f$, but not on the magnitude
of $p_0$.  This principle extends to the case of equations of the form
\begin{equation}
y^{(n)}(t) + p_{n-1}(t) y^{(n-1)}(t) + \cdots + p_1(t) y'(y) + p_0(t) y(t) = f(t).
\label{scalar:inhom}
\end{equation}
That is, such equations admit solutions whose complexity depends on that of the right-hand side $f$ and of
the coefficients $p_0,\ldots,p_{n-1}$, but not on the magnitudes of the coefficients $p_0,\ldots,p_{n-1}$.
This is exploited in \cite{Scalar1} by applying Newton's method to the Riccati equation for (\ref{introduction:scalarode}).
Starting the Newton iterations with a slowly-varying initial guess ensures that each of 
linearized equations have slowly-varying coefficients, and so admit  slowly-varying
solutions.    Consequently, a slowly-varying solution of the Riccati equation can be constructed
via Newton's method as long as an appropriate initial guess is known.
Conveniently enough, there is an obvious mechanism for generating $n$ slowly-varying
initial guesses for the solution of the $(n-1)^{st}$ order Riccati equation.  In particular,
the eigenvalues $\lambda_1(t),\ldots,\lambda_n(t)$ of the coefficient matrix (\ref{introduction:scalarcoef}),
which are often used as low-accuracy approximations of solutions
of the Riccati equation in asymptotic methods,  are suitable as  initial guesses for the Newton procedure.
A proof which parallels this approach is given
in \cite{BremerRokhlin} in order to develop estimates on the complexity of 
the slowly-varying phase functions for  second order equations.

Complicating matters is the fact that the differential operator 
\begin{equation}
D\left[y\right](t) = y^{(n)}(t) + p_{n-1}(t) y^{(n-1)}(t) + \cdots + p_1(t) y'(y) + p_0(t) y(t)
\label{scalar:dop}
\end{equation}
appearing on the left-hand side of   (\ref{scalar:inhom}) 
admits a nontrivial nullspace comprising all solutions of the homogeneous equation
\begin{equation}
y^{(n)}(t) + p_{n-1}(t) y^{(n-1)}(t) + \cdots + p_1(t) y'(y) + p_0(t) y(t) = 0.
\label{scalar:hom}
\end{equation}
This means, of course, that (\ref{scalar:inhom}) is not uniquely solvable.
But it also implies that most solutions of (\ref{scalar:inhom}) are rapidly-varying when the 
coefficients  $p_0,\ldots,p_{n-1}$  are of large magnitude
since the homogeneous equation (\ref{scalar:hom}) admits rapidly-varying solutions in such cases.
It is observed in \cite{Levin} that when the solutions of (\ref{scalar:hom})
are all  rapidly-varying but (\ref{scalar:inhom}) admits a slowly-varying solution $y_0$,
a simple spectral collocation method can be used to compute $y_0$ provided some care is taken in choosing  the discretization
grid.   In particular, if the collocation grid is sufficient to resolve the slowly-varying solution $y_0$, 
but not the solutions of (\ref{scalar:hom}), then the 
matrix discretizing (\ref{scalar:dop}) will be well-conditioned and 
inverting it will yield $y_0$.

The article \cite{Scalar1} describes two algorithms for constructing 
slowly-varying phase functions for scalar linear ordinary differential equations 
based on these principles.  The first such algorithm,
which we refer to as the global Levin method, is an adaptive implementation of the 
procedure described above.  It operates by subdividing $[a,b]$
until, on each resulting subinterval, every one of the slowly-varying phase functions
$\psi_1,\ldots,\psi_n$ it constructs is  accurately represented by a  Chebyshev expansion of a fixed order. The derivatives
$r_1,\ldots,r_n$  of the phase functions $\psi_1,\ldots,\psi_n$
are calculated on each subinterval by applying Newton's method to the Riccati equation and solving the resulting linearized equations via a 
spectral collocation method.   This approach is highly effective as long
as the eigenvalues $\lambda_1(t),\ldots,\lambda_n(t)$ are all of large magnitude on the
interval $[a,b]$.   However, it can fail when one or more of the eigenvalues is of small
magnitude (see \cite{Scalar1} for a further discussion of the problems which arise in this event).

The second algorithm introduced in \cite{Scalar1}, which we refer to as the local Levin method, 
overcomes the difficulties  which arise when one or more of the eigenvalues
of the coefficient matrix (\ref{introduction:scalarcoef})  is of  small magnitude.  
The local method operates in a manner very similar to the algorithm of \cite{BremerPhase}
---  it first computes the values of  the derivatives $r_1, r_2,\ldots, r_n$ of the
desired slowly-varying phase functions $\psi_1,\ldots,\psi_n$  at a point $c$ and then 
solves the Riccati equation numerically with those values used as initial conditions
in order to  construct $r_1,r_2, \ldots,r_n$ over the whole interval $[a,b]$. 
Rather than the continuation method of  \cite{BremerPhase}, however, the values
of the $r_1, r_2,\ldots, r_n$ at $c$ are computed by applying the   Levin approach to a single 
small subinterval of $[a,b]$ containing $c$.     The Riccati equation is then solved
to calculate  $r_1,\ldots,r_n$ over the whole interval, and these
are integrated to obtain the phase functions $\psi_1,\ldots,\psi_n$.
Because most solutions of the Riccati equation are rapidly-varying and we are searching for one of a small
number of slowly-varying solutions, the problems being solved are extremely stiff.
We use a fairly standard adaptive Chebyshev method designed for stiff problems.
It is described in detail in \cite{Scalar1}, but we note that 
essentially any solver for ordinary differential equations should suffice.

The algorithm of this paper uses the local Levin method, which is described in considerably
more detail in \cite{Scalar1}, to construct slowly-varying phase  functions for scalar equations.  
It takes as input
\begin{enumerate}
\item
the interval $[a,b]$ over which the equation is given;

\item
an external subroutine for evaluating the coefficients $q_0,\ldots,q_{n-1}$ in (\ref{introduction:scalarode});

\item
a point $\eta$ on the interval $[a,b]$ and the desired values $\psi_1(\eta),\ldots,\psi_n(\eta)$ for the phase
functions at that point;

\item
an integer $k$ which controls the order of the piecewise Chebyshev expansions used to represent phase functions;

\item
a parameter $\epsilon$ which specifies the desired accuracy for the phase functions; and

\item
a subinterval $[a_0,b_0]$ of $[a,b]$ over which the Levin procedure is to be applied and a point $\sigma$
in that interval.

\end{enumerate}
The output of the algorithm is a collection of  $n^2$ piecewise Chebyshev expansions
of order $(k-1)$  representing the
desired  slowly-varying phase functions $\psi_1,\ldots,\psi_n$ and their derivatives
of orders through $(n-1)$.   By a  $(k-1)^{st}$ order piecewise Chebyshev 
expansions  on the interval $[a,b]$, we mean  a sum of the form
\begin{equation}
\begin{aligned}
&\sum_{i=1}^{m-1} \chi_{\left[x_{i-1},x_{i}\right)} (t) 
\sum_{j=0}^{k-1} \lambda_{ij}\ T_j\left(\frac{2}{x_{i}-x_{i-1}} t + \frac{x_{i}+x_{i-1}}{x_{i}-x_{i-1}}\right)\\
+
&\chi_{\left[x_{m-1},x_{m}\right]} (t) 
\sum_{j=0}^{k-1} \lambda_{mj}\ T_j\left(\frac{2}{x_{m}-x_{m-1}} t + \frac{x_{m}+x_{m-1}}{x_{m}-x_{m-1}}\right),
\end{aligned}
\label{scalar:chebpw}
\end{equation}
where $a = x_0 < x_1 < \cdots < x_m = b$ is a partition of $[a,b]$,
$\chi_I$ is the characteristic function on the interval $I$ and 
$T_j$ is the Chebyshev polynomial of degree $j$.
We note that terms appearing in the first line of (\ref{scalar:chebpw}) involve
the characteristic function of a half-open interval, while that appearing
in the second involves  the characteristic function of a closed interval.
This ensures that exactly one term in  (\ref{scalar:chebpw}) is nonzero for each point $t$ in $[a,b]$.

We note that the matrix
\begin{equation}
\Theta(t) = 
\left(
\begin{array}{ccccc}
u_1(t)  & u_2(t)  & \cdots & u_n(t) \\
u_1'(t) & u_2'(t) & \cdots & u_n'(t) \\
\vdots  & \vdots  & \ddots & \vdots\\
u_1^{(n-1)}(t)  & u_2^{(n-1)}(t)  & \cdots & u_n^{(n-1)}(t) \\
\end{array}
\right),
\end{equation}
where the $u_j(t)$ are given by 
\begin{equation}
u_j(t) = \exp\left(\psi_j(t)\right),
\end{equation}
is a fundamental matrix for the system
\begin{equation}
\mathbf{z}'(t) = B(t) \mathbf{z}(t)
\label{scalar:system1}
\end{equation}
with $B(t)$ is the coefficient matrix (\ref{introduction:scalarcoef})
associated with the scalar equation.  We mean by this that the columns of 
$\Theta(t)$ constitute a basis in the space of solutions of (\ref{scalar:system1}).

\end{section}

\begin{section}{Numerical Algorithm}
\label{section:algorithm}

Here, we describe our numerical algorithm for solving initial and boundary value
problems for a system of linear ordinary differential equations
of the form (\ref{introduction:system}) over an interval $[a,b]$
on which it is nondegenerate.  That is, we  assume that the eigenvalues 
$\lambda_1(t),\ldots,\lambda_n(t)$ of the system's coefficient matrix $A(t)$ are distinct 
on $(a,b)$.  Our algorithm can easily be extended to the case in which the equation
is nondegenerate except at a finite set of turning points in $[a,b]$ by applying
it to a collection of subintervals of $[a,b]$.

Our scheme takes as input the following:
\begin{enumerate}
\item
An integer $n$ specifying the dimension of the system;

\item
the interval $[a,b]$ over which the problem is given;

\item
an integer $k$ which controls the order of the piecewise Chebyshev expansions
used to represent functions on $[a,b]$;

\item
a subroutine for evaluating the elements $a_{ij}(t)$ of the matrix $A(t)$ and
all of their derivatives of orders up to $n$ at a specified point $t$;

\item
a parameter $\epsilon_{\mbox{\tiny disc}}$ specifying the desired precision
for discretizations used to represent the inverse transformation
matrix $\Phi(t)^{-1}$ and the  coefficients $q_0(t),\ldots,q_{n-1}(t)$ 
of the scalar equation;

\item
a parameter $\epsilon_{\mbox{\tiny phase}}$ specifying the desired precision for the
phase functions representing the solutions of the scalar equation;

\item
a subinterval $[a_0,b_0]$ of $[a,b]$ over which the Levin procedure 
used in the construction of the slowly-varying phase functions
for the scalar equation is performed; and 

\item
a constant vector $\mathbf{v} \in \mathbb{C}^n$.

\end{enumerate}
It outputs two collections of piecewise Chebyshev expansions.
The first collection consists of   $n^2$  piecewise Chebyshev expansions of
order $(k-1)$, each of which represents one entry of the inverse  $\Phi(t)^{-1}$ of the transformation matrix 
$\Phi(t)$ defined via (\ref{reduction:transform}) with $\mathbf{v}$ the vector supplied by the user.
The second collection of  expansions comprises $n^2$ piecewise Chebyshev expansions of order $(k-1)$
representing $n$ slowly-varying phase functions $\psi_1(t),\psi_2(t),\ldots,\psi_n(t)$
for the scalar equation 
\begin{equation}
y^{(n)}(t) + q_{n-1}(t) y^{(n-1)}(t) + \cdots + q_1(t) y'(t) + q_0(t) y(t) = 0
\label{algorithm:scalarode}
\end{equation}
corresponding to the coefficient matrix $B(t)$  given by the formula (\ref{reduction:bt}),
as well as  the derivatives up to order $(n-1)$ of these phase functions.  The matrix
\begin{equation}
\Theta(t) = \left(
\begin{array}{cccccc}
u_1(t)          & u_2(t)         & \cdots            &   u_n(t)           & \\
u_1'(t)         & u_2'(t)        & \cdots            &   u_n'(t)            & \\
\vdots          & \vdots         & \ddots      &    \vdots           & \\
u_1^{(n-1)}(t)   & u_2^{(n-1)}(t)   & \cdots            &   u_n^{(n-1)}(t)             & \\
\end{array}
\right),
\label{algorithm:scalarfun}
\end{equation}
where the $u_j$ are given by 
\begin{equation}
u_j(t) = \exp(\psi_j(t)),\ \ \ \ \ j=1,\ldots,n,
\end{equation}
is a fundamental matrix for the scalar system corresponding to (\ref{algorithm:scalarode}), and
\begin{equation}
\Phi(t)^{-1} \Theta(t)
\label{algorithm:fun}
\end{equation}
is a fundamental matrix for the original system (\ref{introduction:system}).
 Once the matrix
$\Phi(t)^{-1} \Theta(t)$ has been formed, essentially any reasonable initial or boundary value
problem for (\ref{introduction:system}) can be solved, essentially instantaneously.


The first step of our algorithm consists of adaptively discretizing $\Phi(t)^{-1}$ and the coefficients 
$q_0(t),\ldots,q_{n-1}(t)$ of the reduced scalar equation.  To do so, our algorithm maintains two lists of 
subintervals of $[a,b]$: one consisting of ``accepted subintervals'' and the other of subintervals which have 
yet to be processed. A subinterval is accepted if the entries of $\Phi(t)^{-1}$ and the functions 
$q_0(t),\ldots,q_{n-1}(t)$ are 
deemed to be adequately represented by a $(k-1)^{st}$ order Chebyshev expansion on that subinterval.
Initially, the list of accepted subintervals is empty and the list of 
subintervals to process contains the single interval $[a,b]$.
We then proceed as follows until the list of subintervals to process is empty:

\begin{enumerate}

\item
Remove a subinterval $[c,d]$ from the list of intervals to process.

\item
Construct the $k$-point extremal Chebyshev grid $t_1,\ldots,t_{k}$ on the
interval $[c,d]$.  The nodes are given by the formula
\begin{equation}
t_j = \frac{d-c}{2} \cos\left(\pi \frac{k-j}{k-1}\right) + \frac{d+c}{2}.
\label{oneint:nodes}
\end{equation}

\item
Evaluate $\Phi(t)$ at each of the nodes $t_1,\ldots,t_k$.  

\item
For each $j=1,\ldots,k$, compute a singular value decomposition
of the matrix $\Phi(t_j)$ and use it to form the inverse $\Phi(t_j)^{-1}$.

\item
For each $j=1,\ldots,k$, use Formula~(\ref{reduction:bt})
to compute the matrix $B(t_j)$.  

\item
Construct $(k-1)^{st}$ order Chebyshev expansions representing each
entry of $\Phi(t)$ and the functions $q_0(t),\ldots,q_{n-1}(t)$,
which appear in the final row of $B(t)$.

\item
For each of the expansions formed in the previous step, which are of the
form
\begin{equation}
\sum_{j=0}^{k-1} \alpha_{j}\ T_j\left(\frac{2}{d-c} t + \frac{d+c}{d-c}\right),
\end{equation}
we compute the ``goodness of fit'' metric
\begin{equation}
\frac
{\sum_{j=k-3}^k\left|\alpha_j\right|^2 }
{\sum_{j=0}^{k-1} \left|\alpha_j\right|^2}.
\end{equation}

\item
If every one of the goodness of fit metrics computed in the previous step is less
than $\epsilon_{\mbox{\tiny disc}}^2$ then  we move the interval $[c,d]$
into the list of accepted intervals.  Otherwise, we put the intervals
\begin{equation}
\left[c,\frac{c+d}{2}\right], \ \ \ \mbox{and} \ \ \left[\frac{c+d}{2},d\right]
\end{equation}
into the list of intervals to process.

\end{enumerate}

Upon termination of this first step, we have piecewise Chebyshev expansions
of order $(k-1)$ representing $\Phi(t)^{-1}$ and the coefficients 
$q_0(t),\ldots,q_{n-1}(t)$ of the reduced scalar equation.  The 
list of accepted intervals determines the partition of the interval $[a,b]$
associated with these piecewise expansions.

We note that the matrix $\Phi(t)$ depends on the entries $a_{ij}(t)$ of the coefficient matrix $A(t)$ and their
derivatives of orders up to $(n-1)$ while $\Phi'(t)$ depends on the $a_{ij}(t)$ and their
derivatives through order $n$.   Because these are supplied by the user as inputs to the algorithm,
there is no need to compute them numerically and suffer the loss of accuracy that repeated numerical
differentiation entails.

In the second step of the algorithm, we use the local Levin method
of \cite{Scalar1}, which is discussed in Section~\ref{section:scalar} of this article,
to  construct piecewise Chebyshev expansions of order $(k-1)$  representing
slowly-varying phase functions $\psi_1(t),\ldots,\psi_n(t)$ for
the scalar equation (\ref{algorithm:scalarode}) 
and their derivatives of orders up to $(n-1)$.  
The precision parameter $\epsilon_{\mbox{\tiny phase}}$ is passed to this algorithm, as is
the integer $k$ controlling the order of the Chebyshev expansions used to represent
function and the interval $[a_0,b_0]$ over which the Levin procedure is to be applied.
For the sake of simplicity, we decided to choose the constants
of integration for the phase function through the requirements that 
\begin{equation}
\psi_1(a) = \psi_2(a) = \cdots = \psi_n(a) = 0.
\end{equation}
The coefficients $q_0(t),\ldots,q_{n-1}(t)$, which are needed for the local Levin method, are
evaluated using the piecewise Chebyshev expansions formed in the first step of the procedure.

Once the second step has been performed, all of the expansions necessary to evaluate
the fundamental matrix for the system (\ref{introduction:system}) have been assembled
and essentially any reasonable initial or boundary value problem for it can be readily solved.

\end{section}

\begin{section}{Numerical Experiments}
\label{section:experiments}

In this section, we present the results of numerical experiments conducted to illustrate the properties
of the algorithms of this paper.  The code for these experiments was written in Fortran and compiled with version 13.1.1
of the GNU Fortran compiler.       They were performed on a desktop computer equipped  with an AMD 3900X processor and 32GB  of memory.  
No attempt  was made to parallelize our code.

The algorithm of \cite{Scalar1}, which we use as a component of the algorithm here,
 calls for computing the eigenvalues of companion matrices.  Standard eigensolvers lose significant 
accuracy when applied to many matrices of this type.   Accordingly, we used the  backward stable and highly accurate method
described in \cite{AURENTZ1,AURENTZ2}  to perform these calculations.

In each of the experiments described below, we considered a system of linear ordinary differential equations
whose coefficient matrix depends on a parameter $\omega$ which controls the frequency of the problem.
For each $\omega = 2^8, 2^9, \ldots, 2^{20}$, we solved an  initial or boundary value problem for the system
over the finite interval $[-1,1]$ using the method of this paper.  
The parameter $k$ controlling the order of Chebyshev expansions used by our algorithm was always taken to be $30$,
while the parameters $\epsilon_{\mbox{\tiny disc}}$ and  $\epsilon_{\mbox{\tiny phase}}$
which specify the desired precision for the discretization of the transformed system and for the solution of the scalar
equation varied from experiment to experiment.  Likewise, different choices of the  interval $[a_0,b_0]$ on which the Levin procedure 
was performed and the  vector $\mathbf{v}$ used to generate  the transformation matrix $\Phi$ were made for each experiment.

The accuracy of each solution $\mathbf{y}$ obtained by our algorithm was measured by comparison with a reference 
solution $\mathbf{z}$ computing using a  standard solver.
To be more explicit, the error in each solution was measured via the quantity 
\begin{equation}
\max_{1 \leq i \leq m} \frac{\left\| \mathbf{y}(t_i) - \mathbf{z}(t_i) \right\|_2}{\left\| \mathbf{z}(t_i) \right\|_2},
\label{experiments:accuracy}
\end{equation}
where $m=10,000$ and $t_1,\ldots,t_{m}$ are the equispaced points in the interval $[-1,1]$ given by the formula
\begin{equation}
t_i = -1 + 2\, \frac{i-1}{m-1}.
\end{equation}
The time taken by our method was measured by repeating the calculation $100$ times and averaging the results.

For each experiment, we provide plots of the time required by our method, the maximum observed value of (\ref{experiments:accuracy})
and the frequency $\Omega$ of the problem as defined in (\ref{introduction:omega}), all as functions of the parameter $\omega$.  
We also report on  one further quantity of interest.  Together, the  phase functions for the transformed scalar equation
and transformation matrix $\Phi^{-1}$ 
represent the basis in the space solutions of (\ref{introduction:system}) used to solve the initial or boundary value
problem.  Each of the entries of $\Phi^{-1}$ and each of the phase functions for the scalar equation
are represented via a piecewise Chebyshev expansion of $(k-1)^{st}$ order.    We provide a plot of the total number of Chebyshev 
coefficients used by all of these expansions, also as a  function of $\omega$.

The conditioning of most initial and boundary value problems for systems of ordinary differential equations
deteriorates as the frequency of the system increases.    Indeed, in typical cases,the 
condition number grows roughly linearly with the frequency $\Omega$ of the system
and  numerical methods show a commensurate decrease in accuracy as the  frequency of the system increases.
So it is not surprising that the experiments here so a similar pattern.

%
%
%
\begin{subsection}{An initial value problem for a system of two equations}
\label{section:experiments:1}

In our first experiment, we solved the system of differential equations
\begin{equation}
\mathbf{y}'(t) = 
\left(
\begin{array}{cc}
1+t^2            & \frac{1}{1+t^4} \\
\frac{-\omega}{1+t^2} & \frac{-i\omega (2+t)}{5+t}
\end{array}
\right)
\mathbf{y}(t)
\label{experiments1:system}
\end{equation}
over the interval $[-1,1]$ subject to the condition
\begin{equation}
\mathbf{y}(0) = \left(
\begin{array}{c}
1\\
1
\end{array}
\right).
\end{equation}
Various computer algebra systems can express the eigenvalues $\lambda_1(t), \lambda_2(t)$ of the coefficient matrix appearing 
in (\ref{experiments1:system}) in terms of elementary functions; however, these formulas are too complicated to reproduce here.  

To give a sense of the magnitudes of $\lambda_1(t), \lambda_2(t)$, we note that 
\begin{equation}
\begin{aligned}
\lambda_1(0) &=  \frac{-2i\omega}{5} - \frac{5i}{2} + \mathcal{O}\left(\frac{1}{\omega}\right)\ \ \  \mbox{and} \\
\lambda_2(0) &=  1 + \frac{5i}{2} + \mathcal{O}\left(\frac{1}{\omega}\right).
\end{aligned}
\end{equation}

The Levin procedure was performed on the interval $[-0.5,0.0]$ and the accuracy parameters 
were taken to be 
\begin{equation}
\epsilon_{\mbox{\tiny disc}} = 1.0 \times 10^{-12} \ \ \ \mbox{and}\ \ \ 
\epsilon_{\mbox{\tiny phase}} = 1.0 \times 10^{-12}.
\end{equation}
The vector $\mathbf{v}$ which generates the transformation matrix $\Phi$ was chosen
to be
\begin{equation}
\mathbf{v} = \left(\begin{array}{c}
1\\
0
\end{array}\right).
\end{equation}

Figure~\ref{experiments1:figure1} gives the results of this experiment, while Figure~\ref{experiments1:figure2}
contains plots of the eigenvalues $\lambda_1(t)$ and $\lambda_2(t)$ when $\omega=2^{16}$.
The frequency $\Omega$ of the systems considered ranged from around 204 (when $\omega=2^8$) to
approximately 821,600 (when $\omega=2^{20}$). 
For every value $\omega$ considered,
exactly 360 Chebyshev coefficients were needed to represent the solutions
of (\ref{experiments1:system}).    
There was surprisingly little variation in the time required by our method
given its complexity and the number of adaptive subprocedures it relies on.  Indeed, approximately  $1.6$ 
milliseconds were required for all values of $\omega$ considered.
As $\omega$ increased from $2^8$ to $2^{20}$, the relative error increased almost linearly
from  $1.5 \times 10^{-13}$ to $4.47 \times 10^{-10}$.

\end{subsection}

%
%
%

\begin{subsection}{A boundary problem for a system of two equations}
\label{section:experiments:2}

In this experiment, we solved the system
\begin{equation}
\mathbf{y}'(t) = 
\left(
\begin{array}{cc}
i \omega\,  \frac{2+\sin^2(6t)}{1+t^2}            & -\frac{\omega}{1+t^2} \\
i + \omega \exp(t)  & i \omega \exp(t)
\end{array}
\right)
\mathbf{y}(t)
\label{experiments2:system}
\end{equation}
over the interval $[-1,1]$ subject to the condition
\begin{equation}
\left(
\begin{array}{cc}
1 & 0 \\
0 & 0
\end{array}
\right)
\mathbf{y}(-1)
+
\left(
\begin{array}{cc}
0 & 0 \\
1 & 0
\end{array}
\right)
\mathbf{y}(1)
=
\left(
\begin{array}{c}
1\\
1
\end{array}
\right).
\label{experiments2:conditions}
\end{equation}
%

\vfil\eject

\begin{figure}[h!]
\hfil
\includegraphics[width=.33\textwidth]{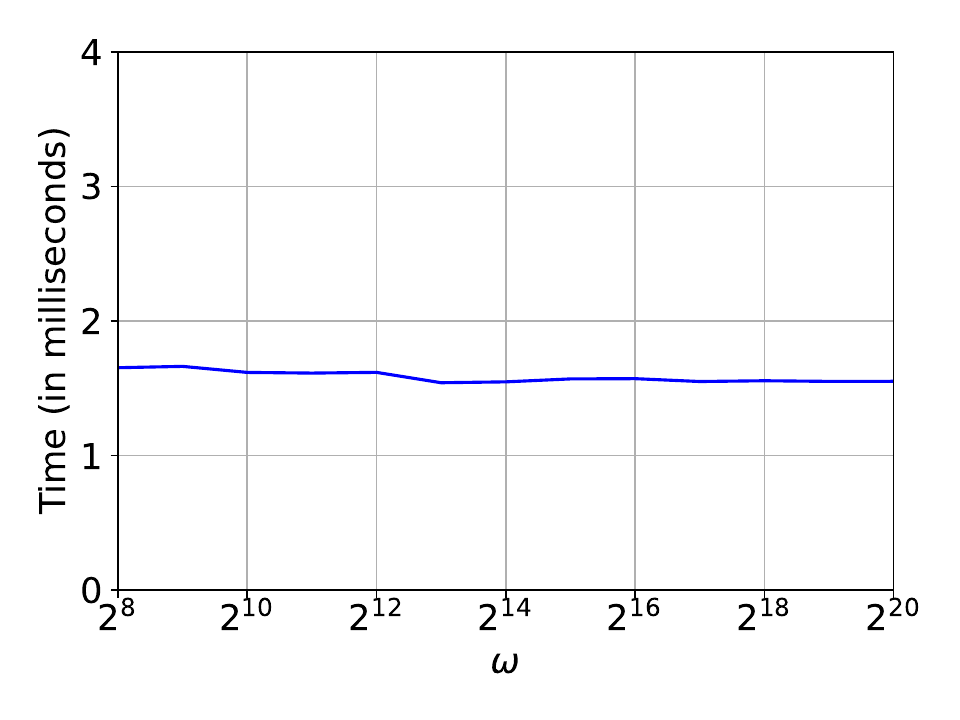}
\hfil
\includegraphics[width=.33\textwidth]{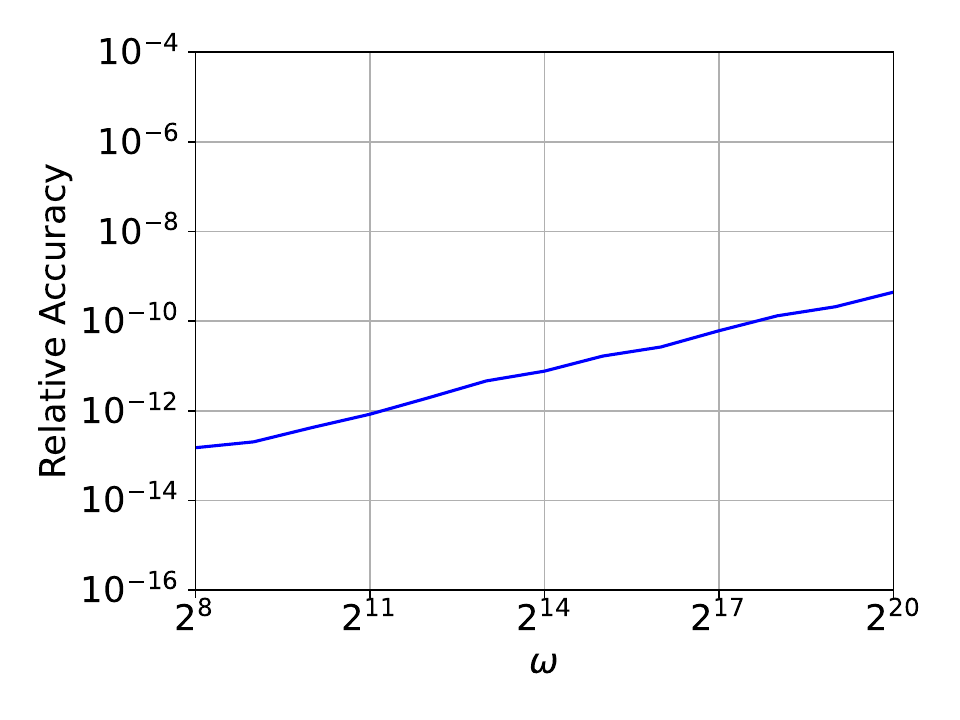}
\hfil

\hfil
\includegraphics[width=.33\textwidth]{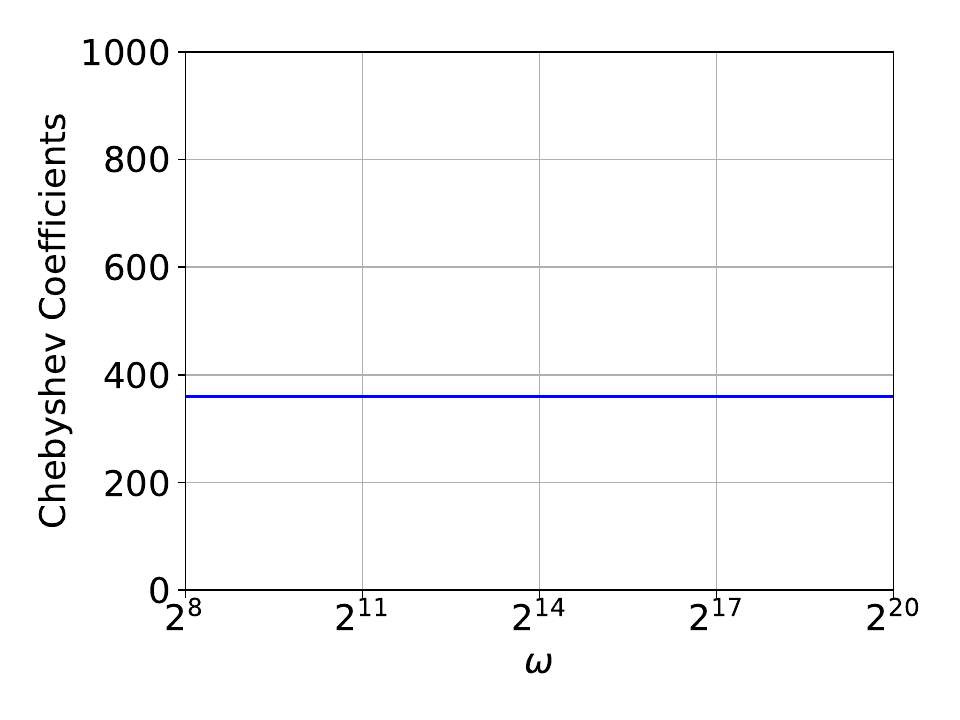}
\hfil
\includegraphics[width=.33\textwidth]{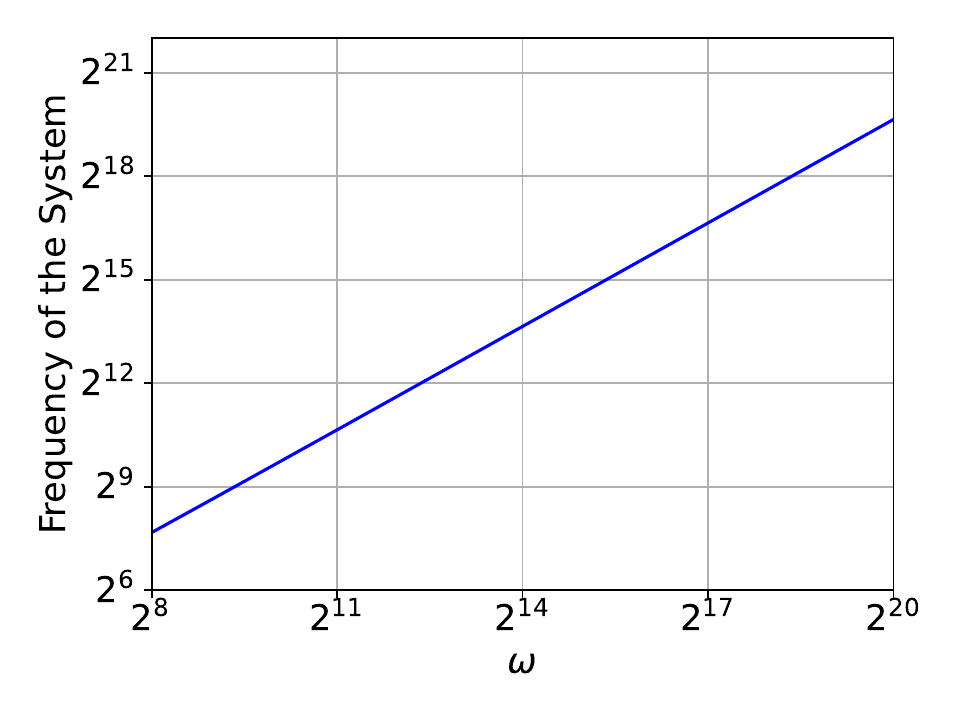}
\hfil

\caption{The results of the experiment of Subsection~\ref{section:experiments:1}.
The plot in the upper left gives the time required by our method
as a function of the parameter $\omega$.  The upper-right plot gives the 
 error in the solution of the initial value problem for the system (\ref{experiments1:system})
as a function of $\omega$.  
The plot on the lower left  shows the total number
of Chebyshev coefficients used to represent the solutions of (\ref{experiments1:system}),
again as a function of the parameter
$\omega$.  The plot on the lower right gives the frequency $\Omega$ of the system as a function of
the parameter $\omega$.}
\label{experiments1:figure1}
\end{figure}
\vfil

\begin{figure}[h!]

\hfil
\includegraphics[width=.33\textwidth]{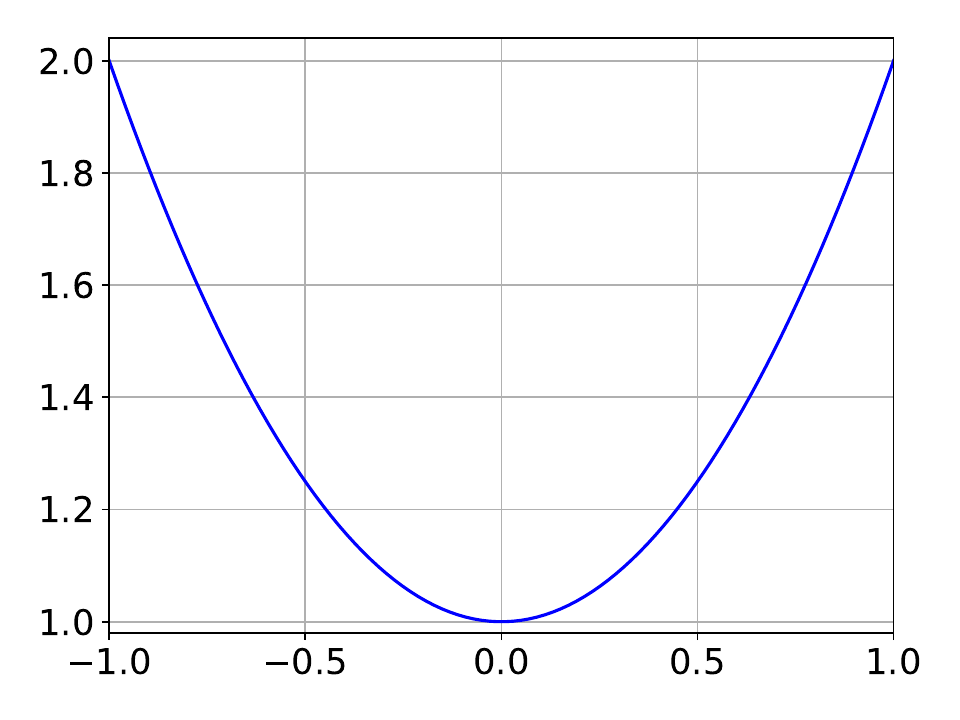}
\hfil
\includegraphics[width=.33\textwidth]{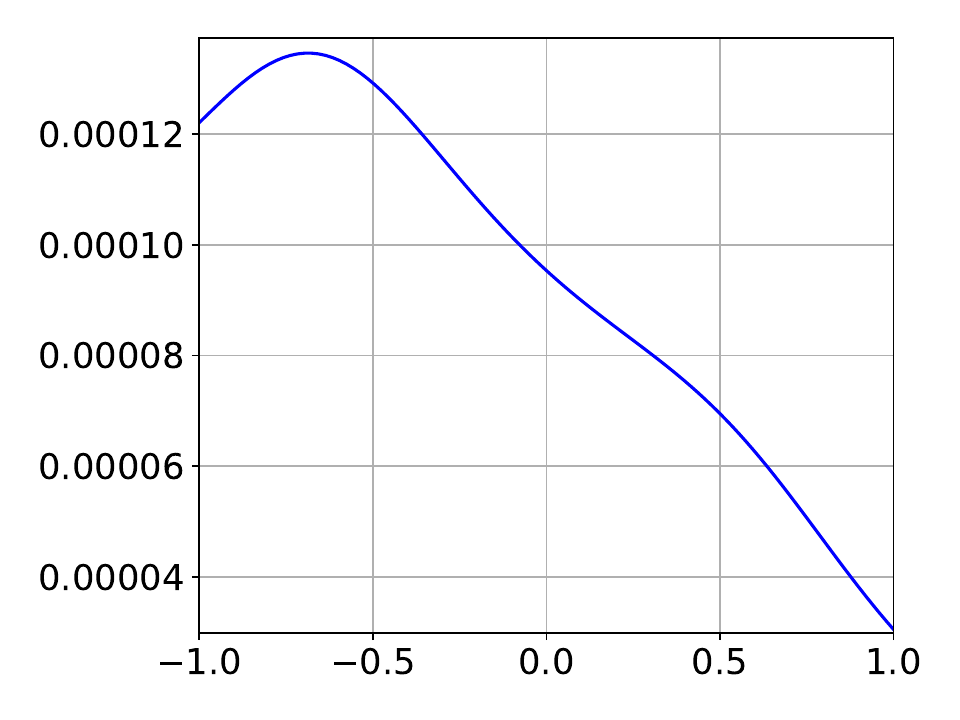}
\hfil

\hfil
\includegraphics[width=.33\textwidth]{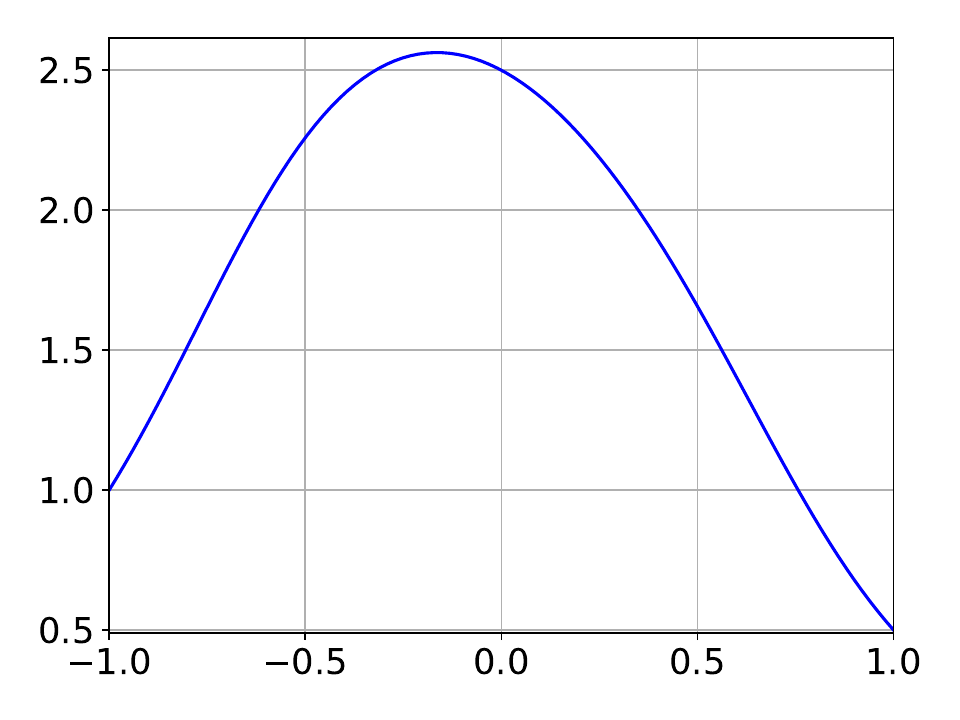}
\hfil
\includegraphics[width=.33\textwidth]{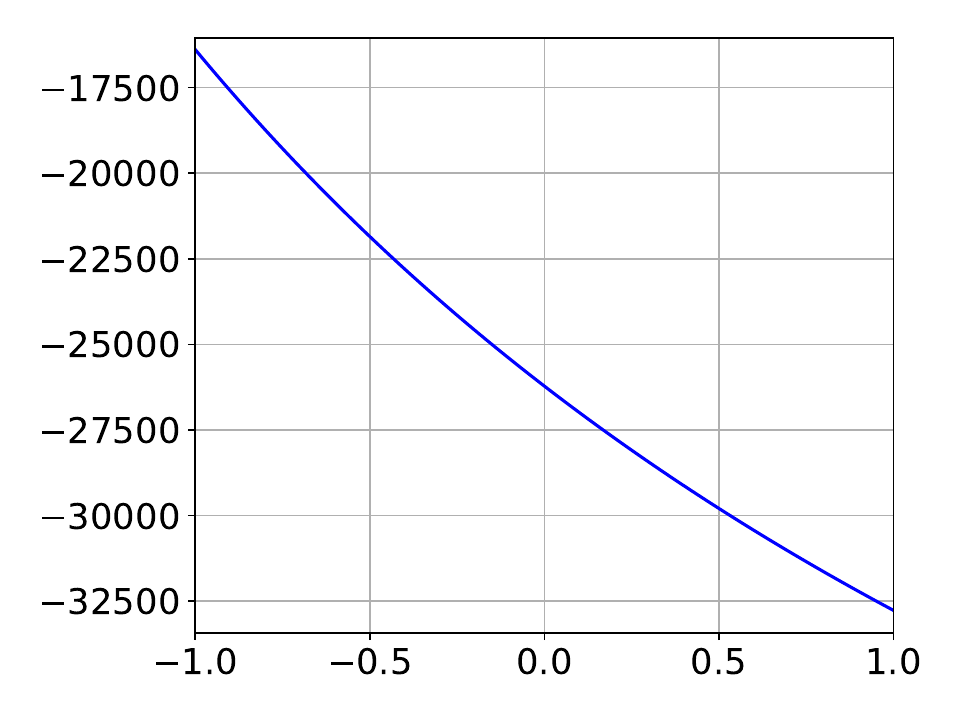}
\hfil

\captionof{figure}{ The eigenvalues $\lambda_1(t)$ and  $\lambda_2(t)$ of the 
coefficient matrix appearing in Equation~(\ref{experiments1:system})
of Subsection~\ref{section:experiments:1}
when the parameter $\omega$ is equal to $2^{16}$.  
Each column corresponds to one eigenvalue,
with the real part appearing in the first row and the imaginary part
in the second. 
}
\label{experiments1:figure2}
\end{figure}

\vfil\eject

Again, while explicit formulas for $\lambda_1(t),\lambda_2(t)$ are available, they are too
complicated to reproduce here.  However, the following formulas give a sense of their
magnitudes:
\begin{equation}
\begin{aligned}
\lambda_1(0) &= \frac{3i}{2} \omega - \frac{i}{2}\sqrt{\omega(4i+5\omega)}\sim 
\frac{i}{2}\left(3+\sqrt{5}\right)\omega - \frac{1}{\sqrt{5}}  + \mathcal{O}\left(\frac{1}{\sqrt{\omega}}\right) \ \ \ \mbox{and}\\
\lambda_2(0) &= \frac{3i}{2} \omega + \frac{i}{2}\sqrt{\omega(4i+5\omega)}\sim
\frac{i}{2}\left(3- \sqrt{5}\right)\omega + \frac{1}{\sqrt{5}}
+ \mathcal{O}\left(\frac{1}{\sqrt{\omega}}\right).
\end{aligned}
\end{equation}

The Levin procedure was performed on the interval $[-0.5,0.0]$
and we took the accuracy parameters to be 
\begin{equation}
\epsilon_{\mbox{\tiny disc}} = 1.0 \times 10^{-12} \ \ \ \mbox{and}\ \ \ 
\epsilon_{\mbox{\tiny phase}} = 1.0 \times 10^{-12}.
\end{equation}
The vector $\mathbf{v}$ which generates the transformation matrix $\Phi$ was
\begin{equation}
\mathbf{v} = \left(\begin{array}{c}
0\\
1
\end{array}\right).
\end{equation}

Figure~\ref{experiments2:figure1} gives the results of this experiment, while Figure~\ref{experiments2:figure2}
contains plots of the eigenvalues $\lambda_1(t)$ and $\lambda_2(t)$ when $\omega=2^{16}$.
The frequency $\Omega$ of the systems considered ranged from around 1,370 (when $\omega=2^8$) to
approximately 5,612,000 (when $\omega=2^{20}$).  When $\omega=2^8$,
900 Chebyshev coefficients were used to represent the solutions
of (\ref{experiments2:system}), while 810 coefficients were needed when $\omega=2^9$
and, for all other values of $\omega$, 720 Chebyshev coefficients sufficed.

This is unsurprising given that the complexity of phase functions generally decreases
with increasing frequency.  The time required by our algorithm was a bit higher 
when $\omega=2^8$ and $\omega=2^9$, but, again, this is a consequence
of the slightly higher complexity of the phase functions for those
values of $\omega$.  

\end{subsection}

%
%
%
\begin{subsection}{An initial value problem for a system of three equations}
\label{section:experiments:3}

In this experiment, we solved the system of differential equations
\begin{equation}
\mathbf{y}'(t) = A(t) \mathbf{y}(t)
\label{experiments3:system}
\end{equation}
where
\begin{equation*}
A(t) = 
\left(
\begin{array}{cc}
 -i k \left(3 t^2+e^{t-12 t^2}+3 e^t+e^t \cos (17 t)+3\right) & -i k \left(e^{-12t^2}+\cos (17 t)+3\right) \\
 -i k e^t \left(-3 t^2+e^{-12 t^2}-2\right) & -i k \left(e^{-12 t^2}+1\right) \\
 -i k e^t \left(e^{-12 t^2}+\cos (17 t)+3\right) & -i k \left(e^{-12 t^2}+\cos (17t)+3\right) \\
\end{array}
\right)
\end{equation*}
\begin{equation*}
\left.
\begin{array}{c}
i k \left(3 t^2+e^{t-12 t^2}+3 e^t+\left(e^t+1\right)\cos (17 t)+5\right) \\
 i ke^t \left(-3 t^2+e^{-12 t^2}-2\right) \\
 i k \left(e^{t-12 t^2}+3 e^t+\left(e^t+1\right) \cos (17t)+2\right) 
\end{array}
\right),
\end{equation*}
\vfil\eject

\begin{figure}[h!]
\hfil
\includegraphics[width=.33\textwidth]{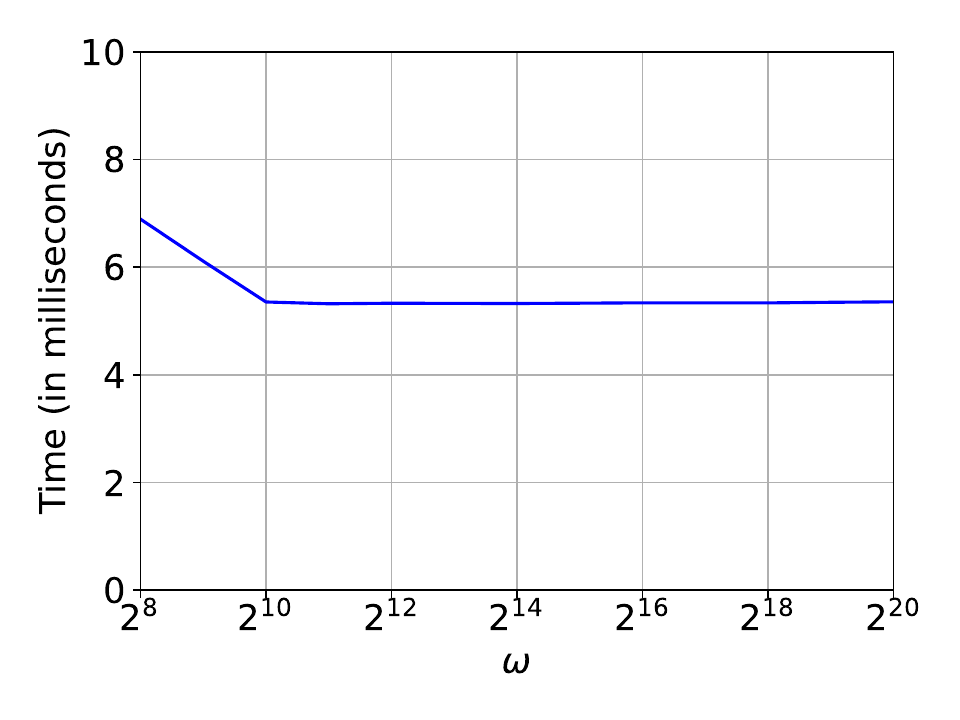}
\hfil
\includegraphics[width=.33\textwidth]{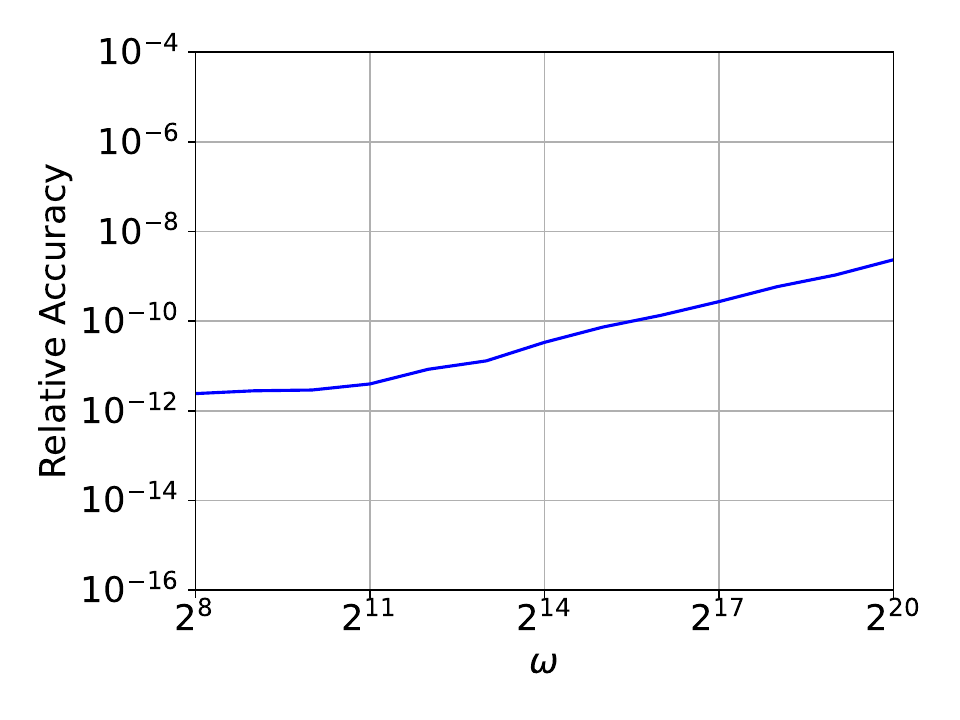}
\hfil

\hfil
\includegraphics[width=.33\textwidth]{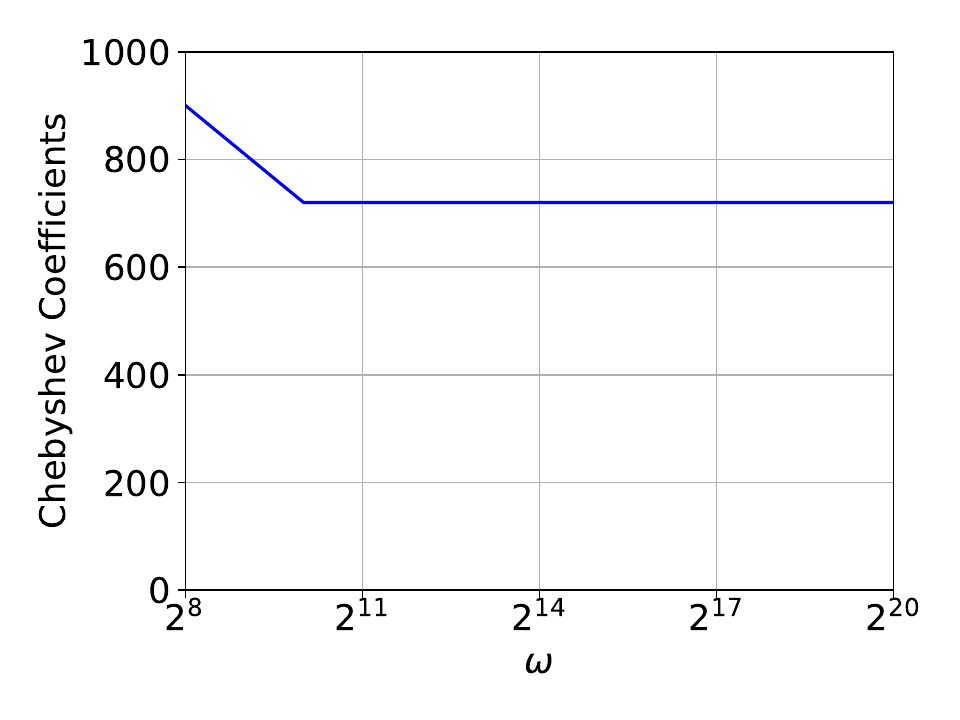}
\hfil
\includegraphics[width=.33\textwidth]{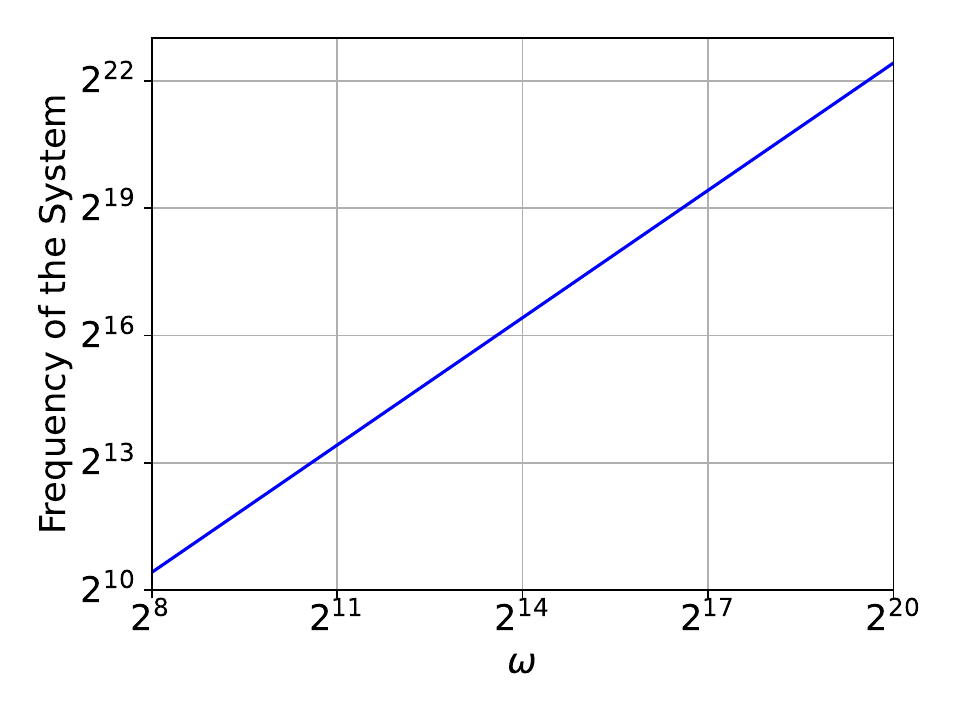}
\hfil

\caption{The results of the experiment of Subsection~\ref{section:experiments:2}.
The plot in the upper left gives the time required by our method
as a function of the parameter $\omega$.  The upper-right plot gives the 
 error in the solution of the initial value problem for the system (\ref{experiments2:system})
as a function of $\omega$.  
The plot on the lower left  shows the total number
of Chebyshev coefficients used to represent the solutions of (\ref{experiments2:system}),
again as a function of the parameter $\omega$.
formed during the experiments  of Section~\ref{section:experiments:2}, again as a function of the parameter
$\omega$.  The plot on the lower right gives the frequency $\Omega$ of the system as a function of
the parameter $\omega$.}
\label{experiments2:figure1}
\end{figure}

\vfil

\begin{figure}[h!]

\hfil
\includegraphics[width=.33\textwidth]{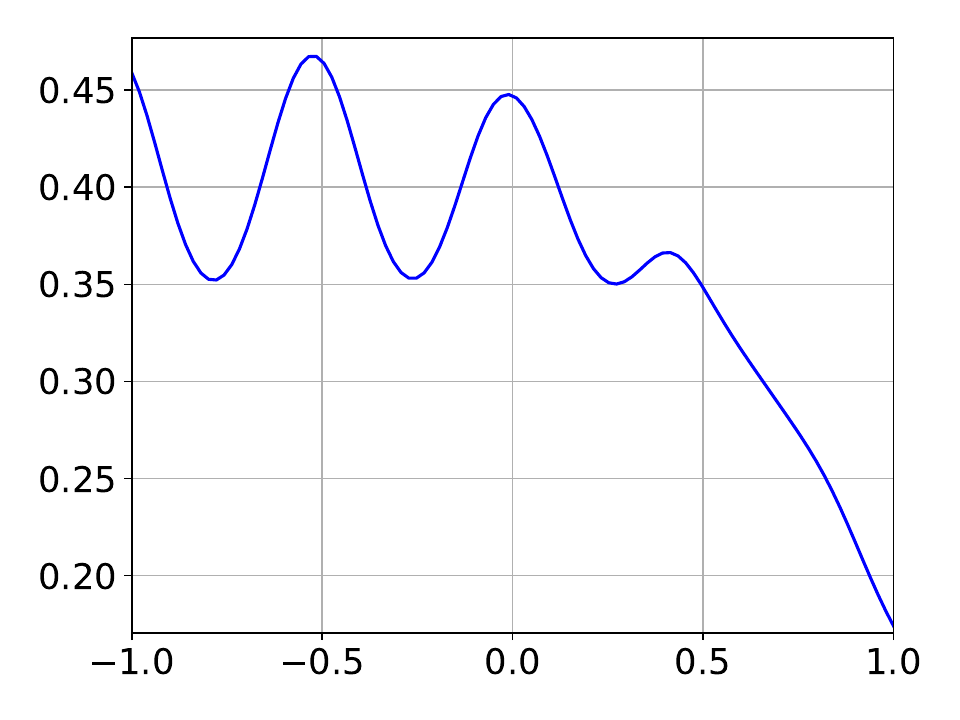}
\hfil
\includegraphics[width=.33\textwidth]{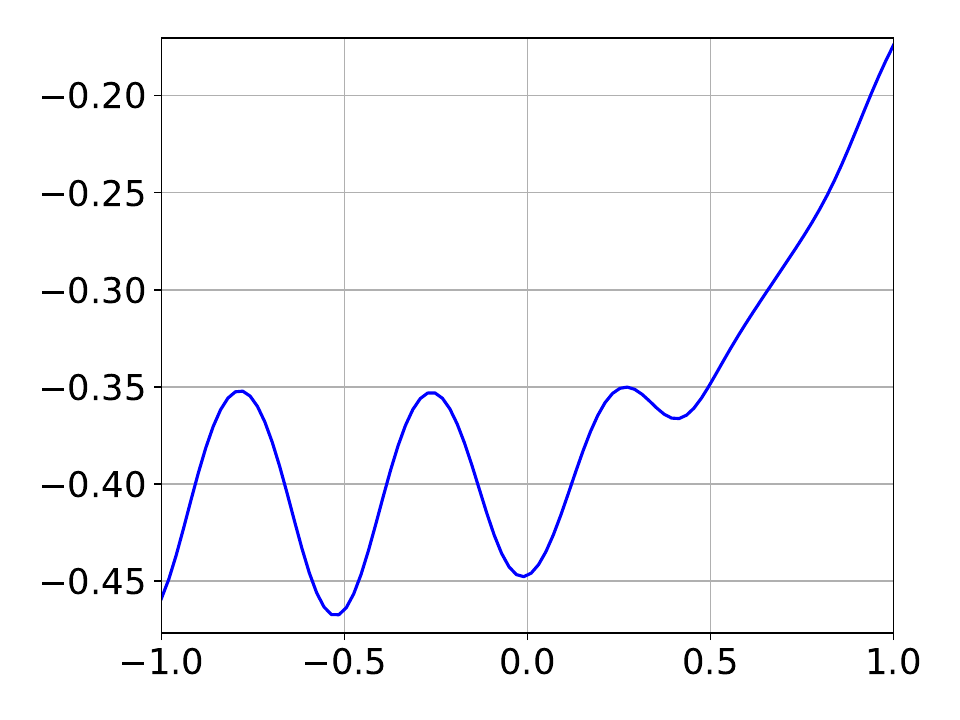}
\hfil

\hfil
\includegraphics[width=.33\textwidth]{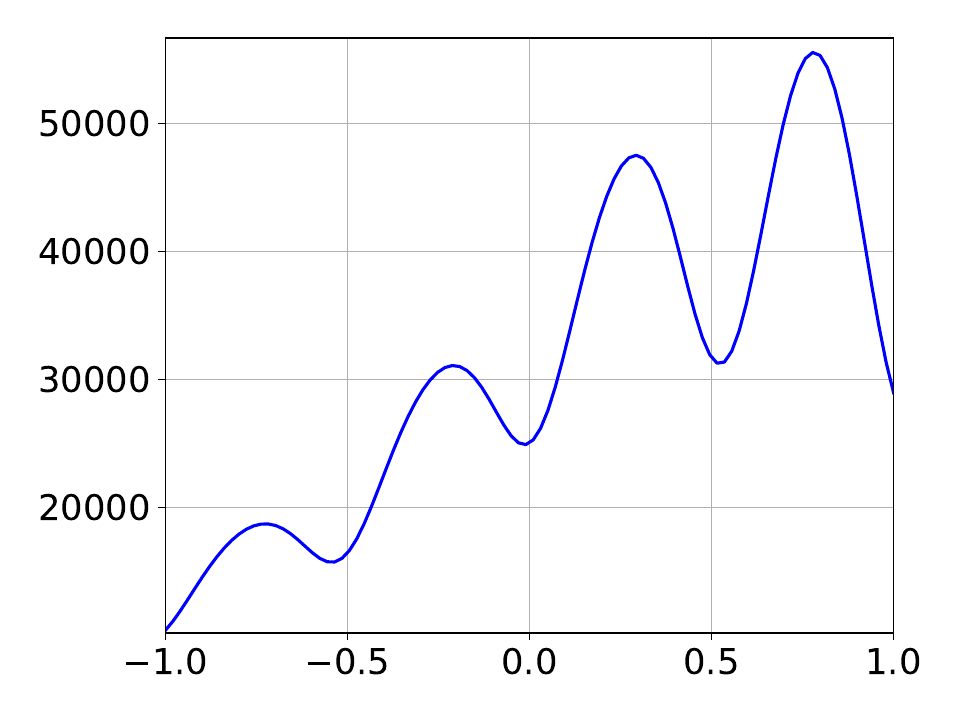}
\hfil
\includegraphics[width=.33\textwidth]{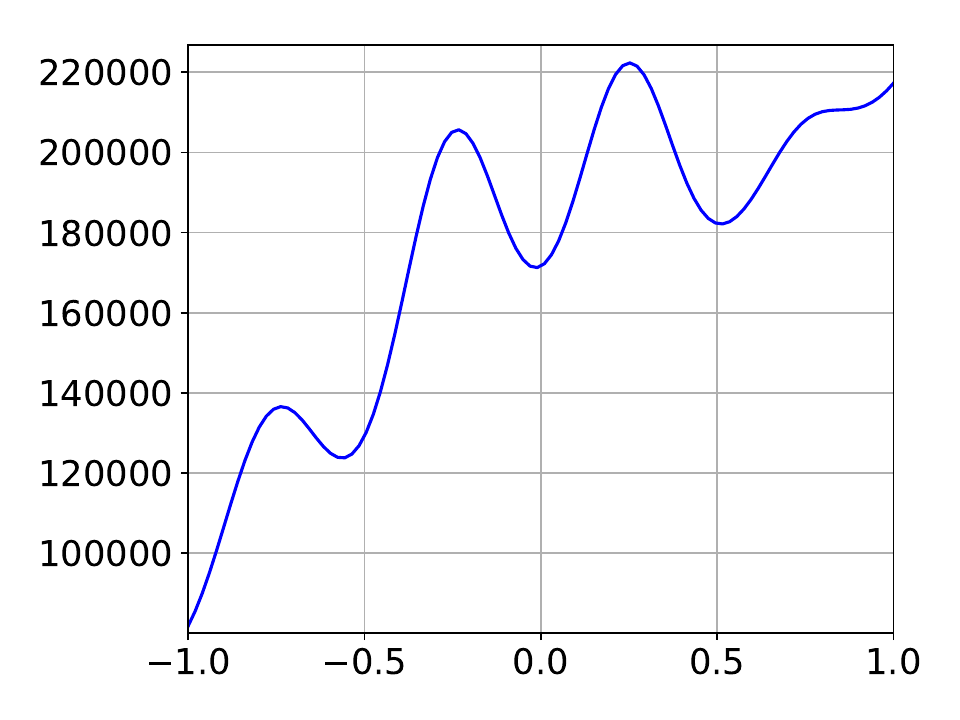}
\hfil

\captionof{figure}{The eigenvalues $\lambda_1(t), \lambda_2(t)$ of the 
coefficient matrix appearing in Equation~(\ref{experiments2:system})
of Subsection~\ref{section:experiments:2}
when the parameter $\omega$ is equal to $2^{16}$.  
Each column corresponds to one eigenvalue,
with the real part appearing in the first row and the imaginary part
in the second. 
}
\label{experiments2:figure2}

\end{figure}

\vfil\eject

over the interval $[-1,1]$ subject to the condition
\begin{equation}
\mathbf{y}(-1) = \left(\begin{array}{c}
1\\
0\\
-1
\end{array}\right).
\end{equation}

The eigenvalues of $A(t)$ are given by the following formulas:
\begin{equation}
\begin{aligned}
\lambda_1(t) &= i\omega \left(2+\cos(12t)\right),\\
\lambda_2(t) &= -3i \omega (1+t^2) \ \ \ \mbox{and} \\
\lambda_3(t) &= -i\omega \left(1+\exp\left(-12t^2\right)\right).
\end{aligned}
\end{equation}

The Levin procedure was performed on the interval $[-0.25,0.0]$ and the accuracy parameters 
were taken to be 
\begin{equation}
\epsilon_{\mbox{\tiny disc}} = 1.0 \times 10^{-12} \ \ \ \mbox{and}\ \ \ 
\epsilon_{\mbox{\tiny phase}} = 1.0 \times 10^{-12}.
\end{equation}
The initial vector for constructing the transformation matrix was chosen to be
\begin{equation}
\mathbf{v} = \left(
\begin{array}{c}
1\\
0\\
0
\end{array}
\right).
\end{equation}
Figure~\ref{experiments3:figure1} gives the results of this experiment, while Figure~\ref{experiments3:figure2}
contains plots of the eigenvalues $\lambda_1(t)$, $\lambda_2(t)$ and $\lambda_3(t)$ when $\omega=2^{16}$.
The frequency $\Omega$ of the systems considered varied from approximately
2,048  (when $\omega=2^8$) to around 8,388,600 (when $\omega=2^{20}$).
When $\omega=2^8$ and $\omega=2^9$, 4,080 Chebyshev coefficients were required to represent the solutions
of (\ref{experiments3:system}), while 4,020 coefficients were needed for values of  $\omega$ between $2^{10}$
and $2^{13}$.  For all other values of $\omega$, 3,960 Chebyshev coefficients were required.
 The time taken by our algorithm also decreased slightly with increasing $\omega$.
There was somewhat more variability in the obtained accuracy than in previous experiments; however, 
relatively high accuracy was still achieved in all cases and it is unclear if the variation is due to 
intrinsic properties of the problem (e.g., variation in the condition number)
or the properties of our algorithm.

\end{subsection}

\vfil\eject

\begin{figure}[h!]
\hfil
\includegraphics[width=.33\textwidth]{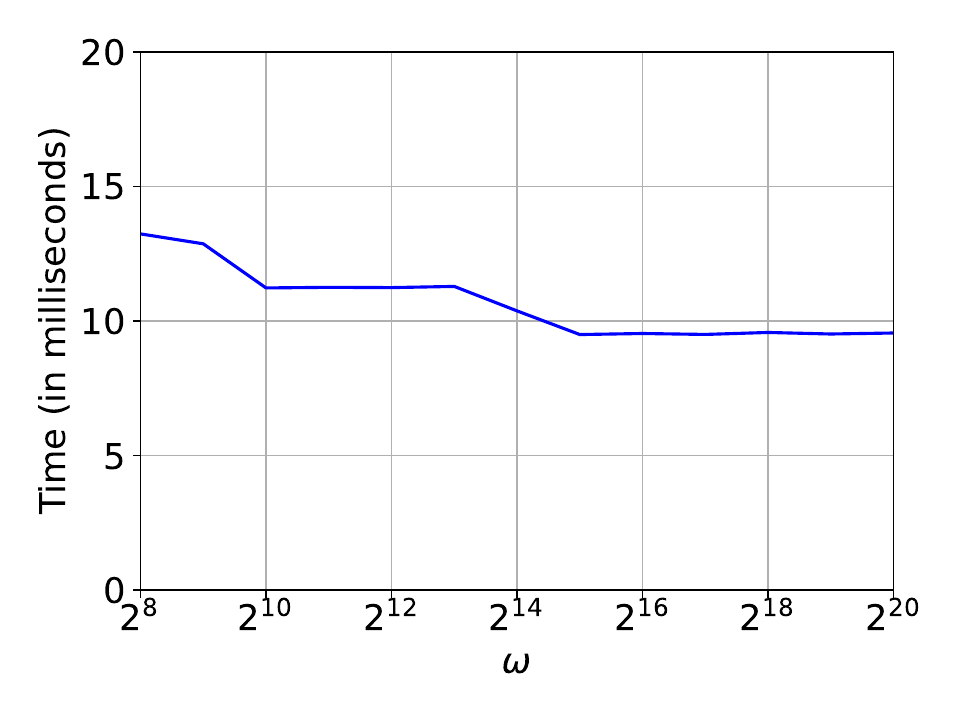}
\hfil
\includegraphics[width=.33\textwidth]{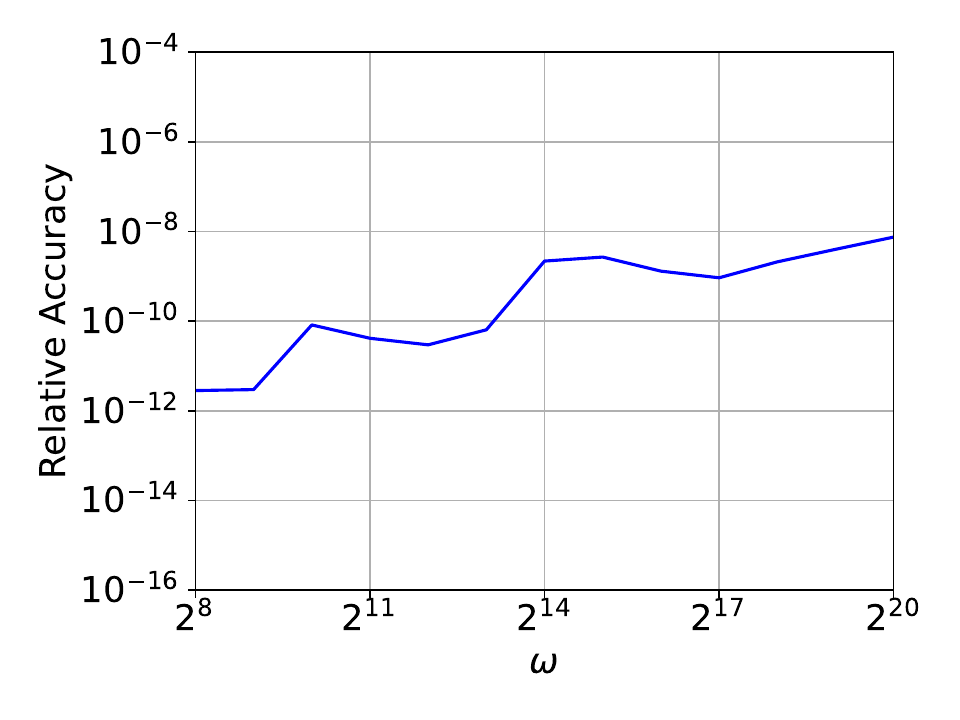}
\hfil

\hfil
\includegraphics[width=.33\textwidth]{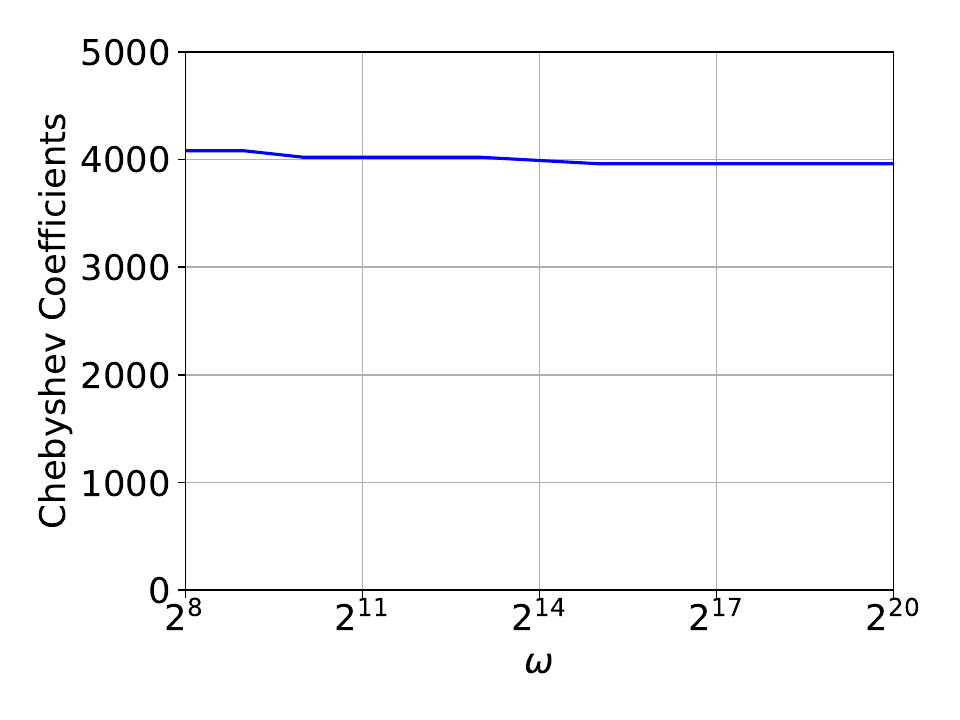}
\hfil
\includegraphics[width=.33\textwidth]{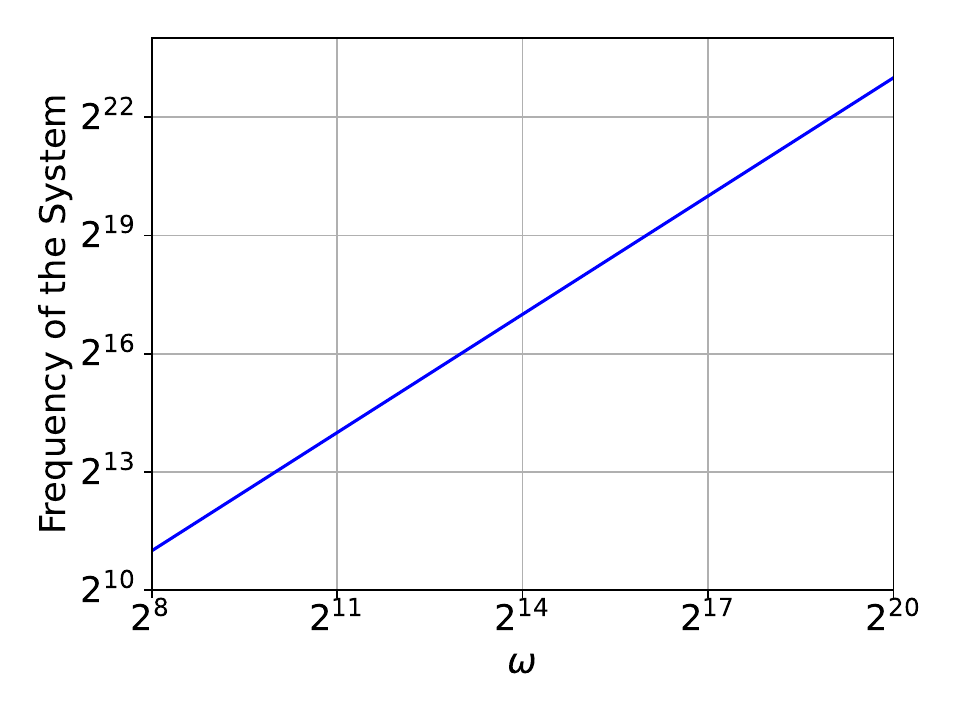}
\hfil

\caption{The results of the experiment of Subsection~\ref{section:experiments:3}.
The plot in the upper left gives the time required by our method
as a function of the parameter $\omega$.  The upper-right plot gives the 
 error in the solution of the initial value problem for the system (\ref{experiments3:system})
as a function of $\omega$.  
The plot on the lower left  shows the total number
of Chebyshev coefficients used to represent the solutions of (\ref{experiments3:system}),
again as a function of the parameter $\omega$.
 The plot on the lower right gives the frequency $\Omega$ of the system as a function of
the parameter $\omega$.}
\label{experiments3:figure1}
\end{figure}

\vfil

\begin{figure}[h!]
\centering
\hfil
\includegraphics[width=.32\textwidth]{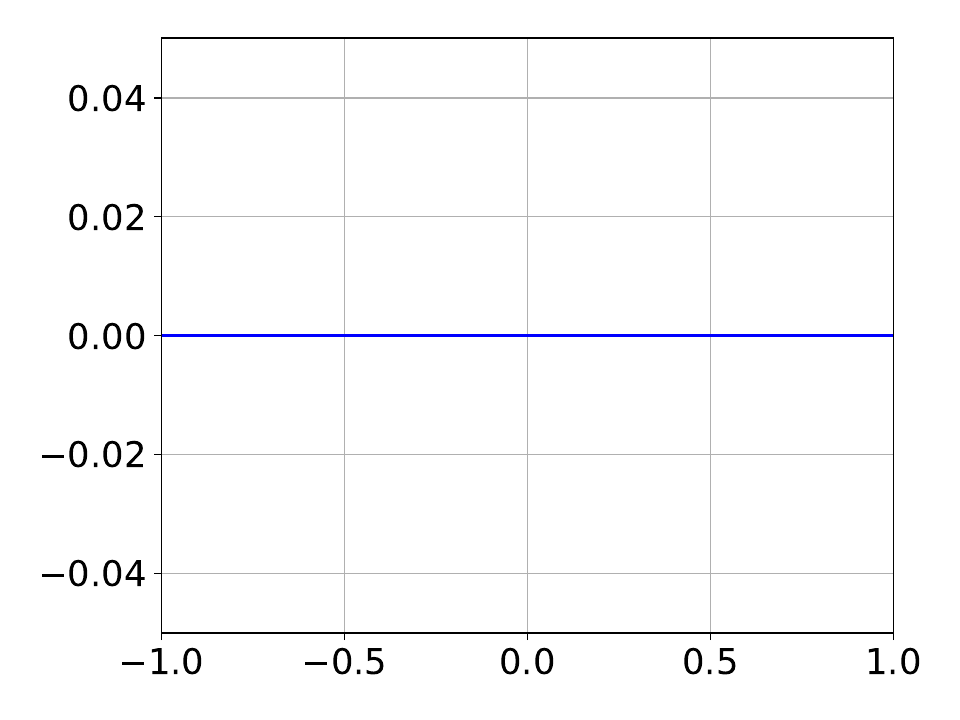}
\hfil
\includegraphics[width=.32\textwidth]{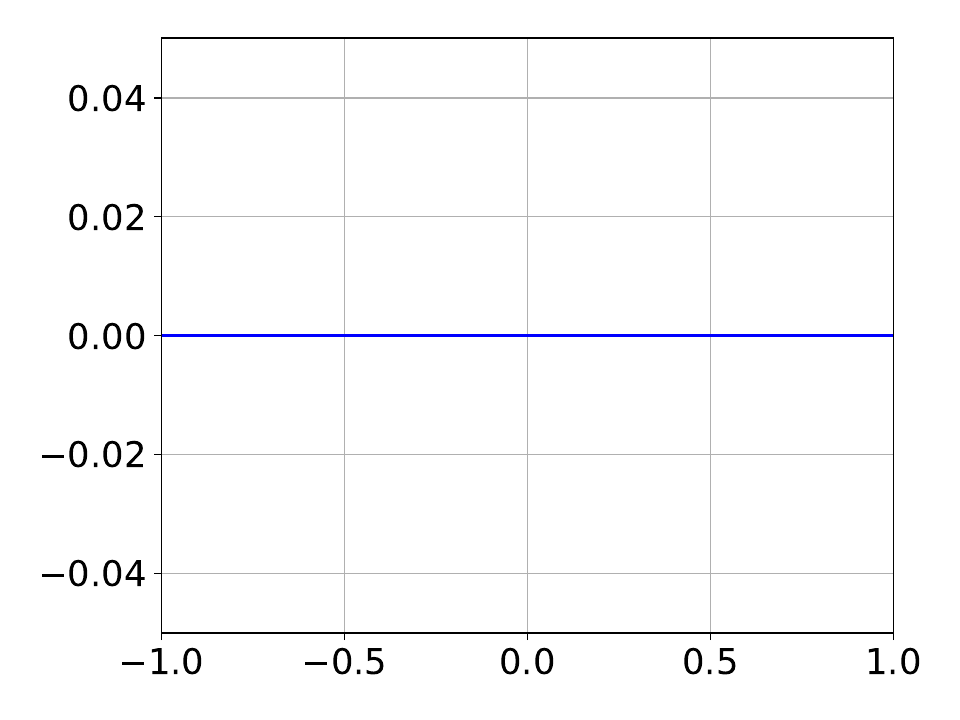}
\hfil
\includegraphics[width=.32\textwidth]{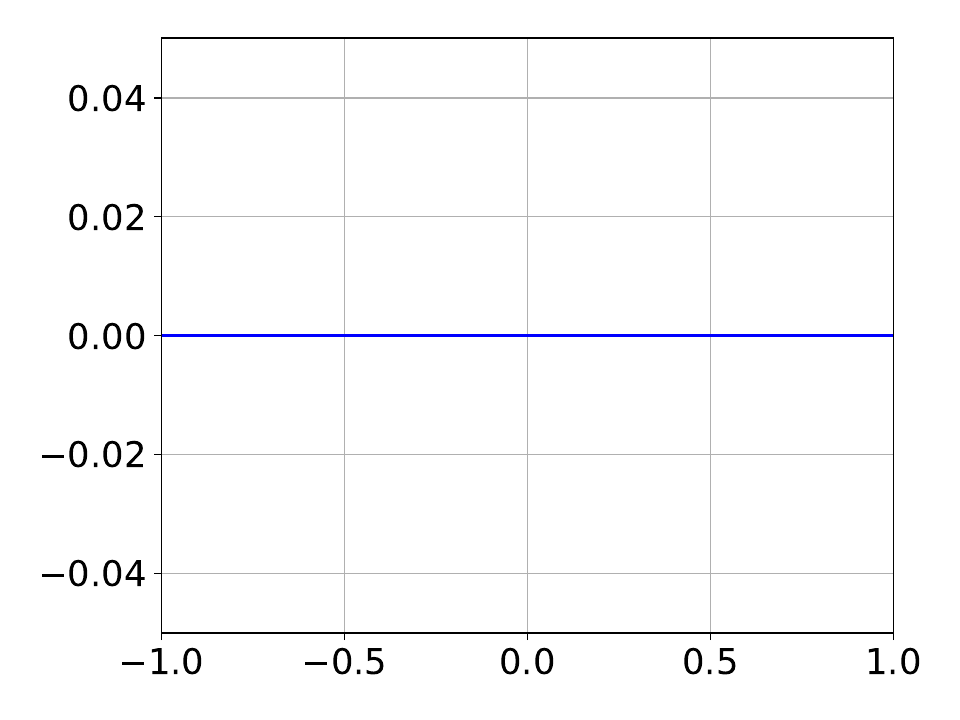}
\hfil

\hfil
\includegraphics[width=.32\textwidth]{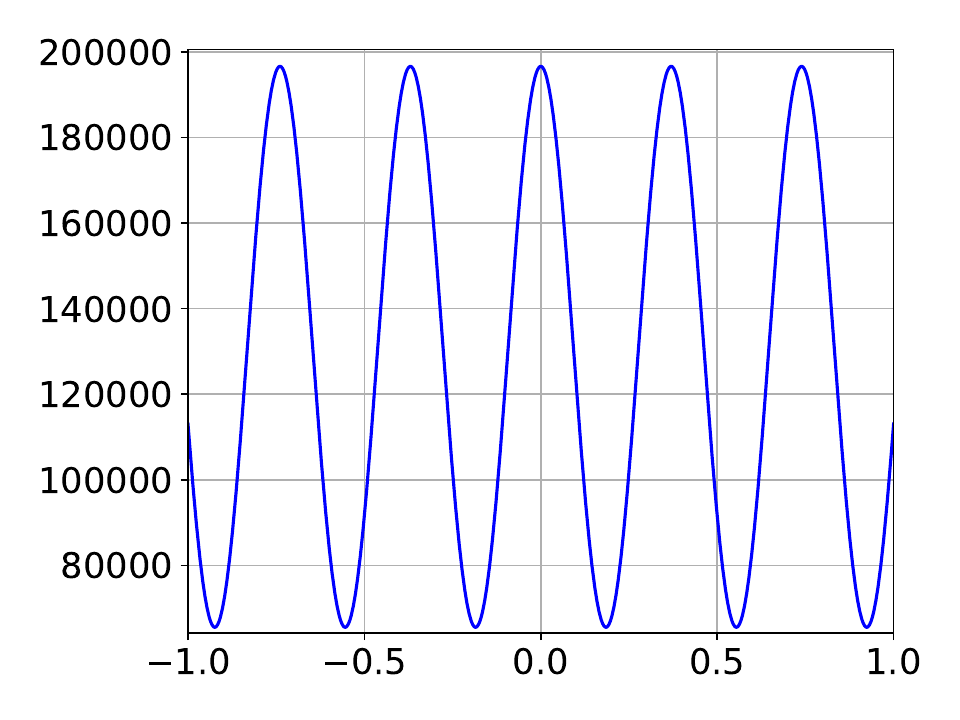}
\hfil
\includegraphics[width=.32\textwidth]{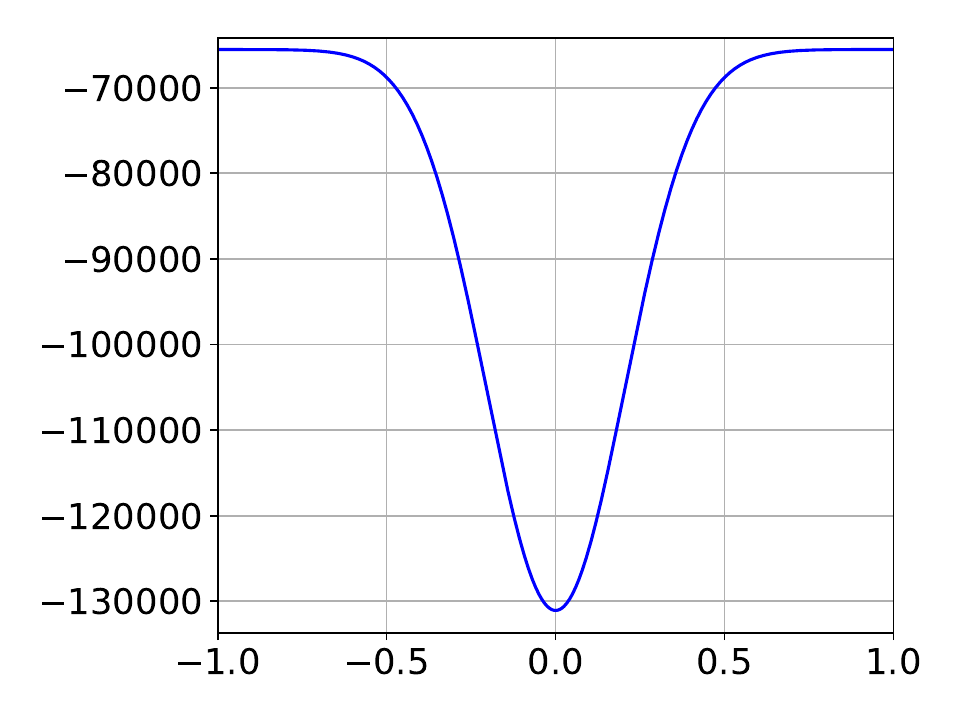}
\hfil
\includegraphics[width=.32\textwidth]{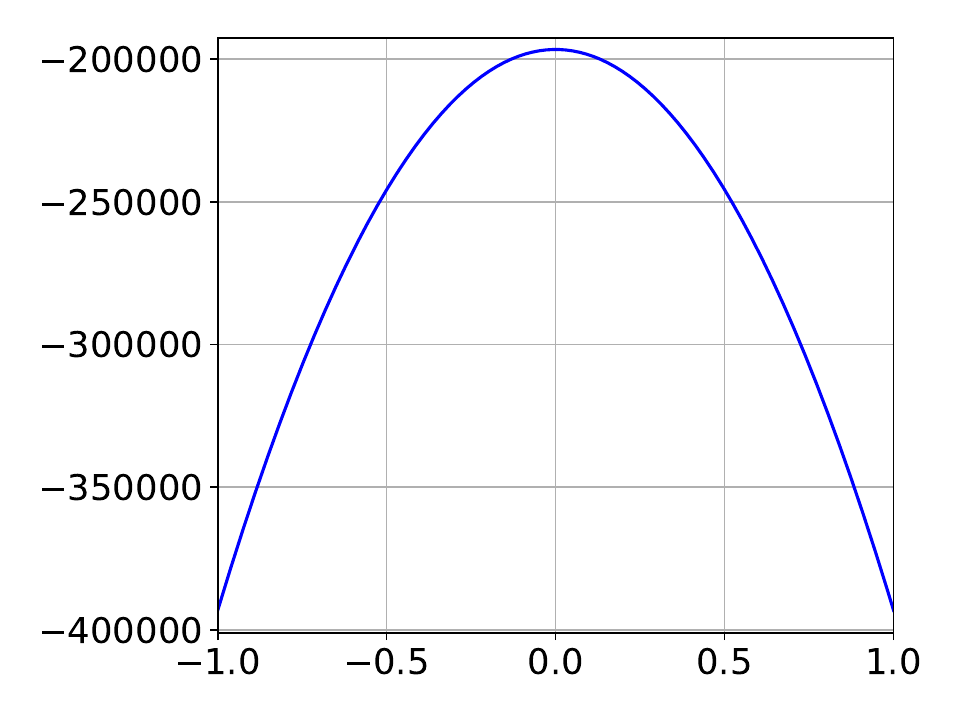}
\hfil

\captionof{figure}{The eigenvalues $\lambda_1(t), \lambda_2(t), \lambda_3(t)$ of the 
coefficient matrix  $A(t)$
of  Subsection~\ref{section:experiments:3}
when the parameter $\omega$ is equal to $2^{16}$.
Each column corresponds to one eigenvalue,
with the real part appearing in the first row and the imaginary part
in the second. 
}
\label{experiments3:figure2}

\end{figure}
\vfil\eject

%
%
%
\begin{subsection}{A boundary value problem for a system of three equations}
\label{section:experiments:4}

In this experiment, we solved the system of differential equations
\begin{equation}
\mathbf{y}'(t) = A(t) \mathbf{y}(t),
\label{experiments4:system}
\end{equation}
where 
\begin{equation*}
A(t) = \left(
\begin{array}{ccc}
 \frac{4 (1+4 i k) e^{3 t^2}+8 i k e^{3 t^2} \sin (3 t)-i k t \log
   \left(t+\frac{1001}{1000}\right)}{4 e^{3 t^2}-t} & 0 \\
 \frac{2 e^{2 t^2} t \left(-8 i k e^{t^2}+2 i k \sin (3 t)+\log (k) \sin (t)+4 i
   k+1\right)}{4 e^{3 t^2}-t} & -\log (k) \sin (t)+8 i k e^{t^2} 
   \\
 \frac{2 e^{2 t^2} \left(-i k \log \left(t+\frac{1001}{1000}\right)+2 i k \sin (3
   t)+4 i k+1\right)}{4 e^{3 t^2}-t} & 0 \\
\end{array}
\right.
\end{equation*}
\begin{equation}
\left.
\begin{array}{c}
 \frac{2 i e^{t^2} t
   \left(k \log \left(t+\frac{1001}{1000}\right)-2 k \sin (3 t)-4 k+i\right)}{4 e^{3t^2}-t} \\
\frac{i t^2
   \left(8 k e^{t^2}-2 k \sin (3 t)+i \log (k) \sin (t)-4 k+i\right)}{4 e^{3 t^2}-t}
   \\
 \frac{i \left(4 k e^{3 t^2} \log
   \left(t+\frac{1001}{1000}\right)+(-4 k+i) t-2 k t \sin (3 t)\right)}{4 e^{3t^2}-t} \\
\end{array}
\right),
\end{equation}
over the interval $[-1,1]$ subject to the condition
\begin{equation}
\left(
\begin{array}{ccc}
1 & 1 & 0\\
1 & 0 & 1\\
0 & 1 & 0\\
\end{array}
\right)
\mathbf{y}(-1)
+
\left(
\begin{array}{ccc}
0 & 0 & 1\\
0 & 1 & 0\\
0 & -1 & 0\\
\end{array}
\right)
\mathbf{y}(1)
=
\left(
\begin{array}{c}
1\\
0\\
1
\end{array}
\right).
\label{experiments4:conditions}
\end{equation}
The eigenvalues of $A(t)$ are given by the following formulas:
\begin{equation}
\begin{aligned}
\lambda_1(t) &= 1+2i\omega\left(2+\sin(3t)\right),\\
\lambda_2(t) &= -\log(k)\sin(t) + 8ik \exp(t^2)\ \ \ \mbox{and}\\
\lambda_3(t) &= ik\log\left(1+10^{-3}+t\right).
\end{aligned}
\end{equation}

The Levin procedure was performed on the interval $[-0.1,0.0]$ and the accuracy parameters 
were taken to be 
\begin{equation}
\epsilon_{\mbox{\tiny disc}} = 1.0 \times 10^{-12} \ \ \ \mbox{and}\ \ \ 
\epsilon_{\mbox{\tiny phase}} = 1.0 \times 10^{-12}.
\end{equation}
As with all of the other experiments discussed here, the parameter $k$ controlling
the size of the piecewise Chebyshev expansions used to represent functions was set to $30$.
We used  
\begin{equation}
\mathbf{v} = \left(
\begin{array}{c}
1\\
1\\
1
\end{array}
\right).
\end{equation}
as the initial vector for constructing the transformation matrix $\Phi$.

\vfil
\eject

\begin{figure}[h!]
\hfil
\includegraphics[width=.33\textwidth]{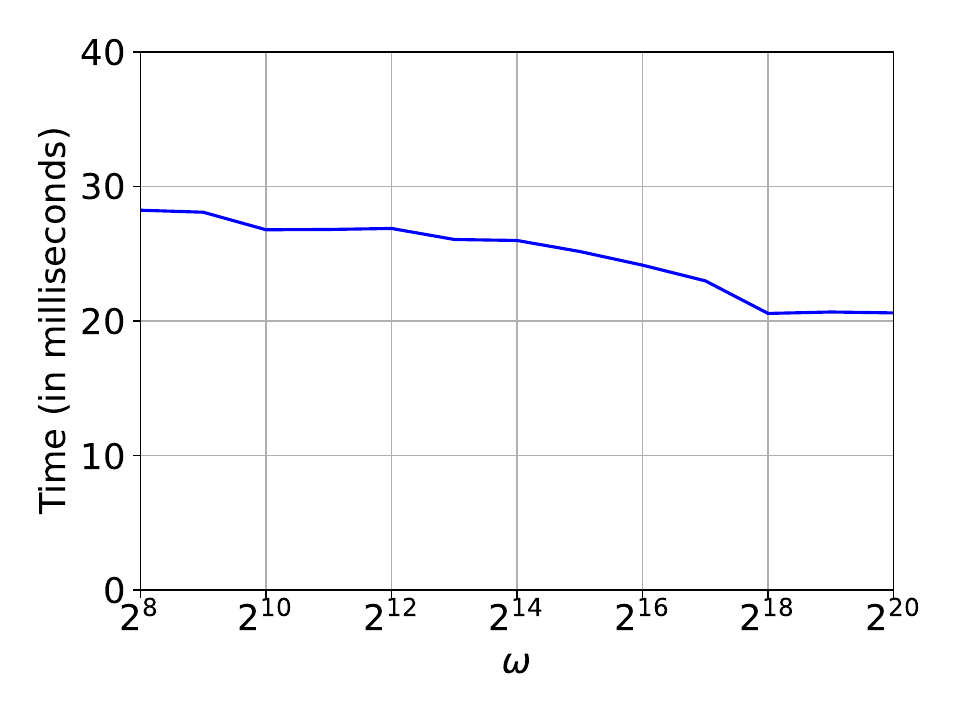}
\hfil
\includegraphics[width=.33\textwidth]{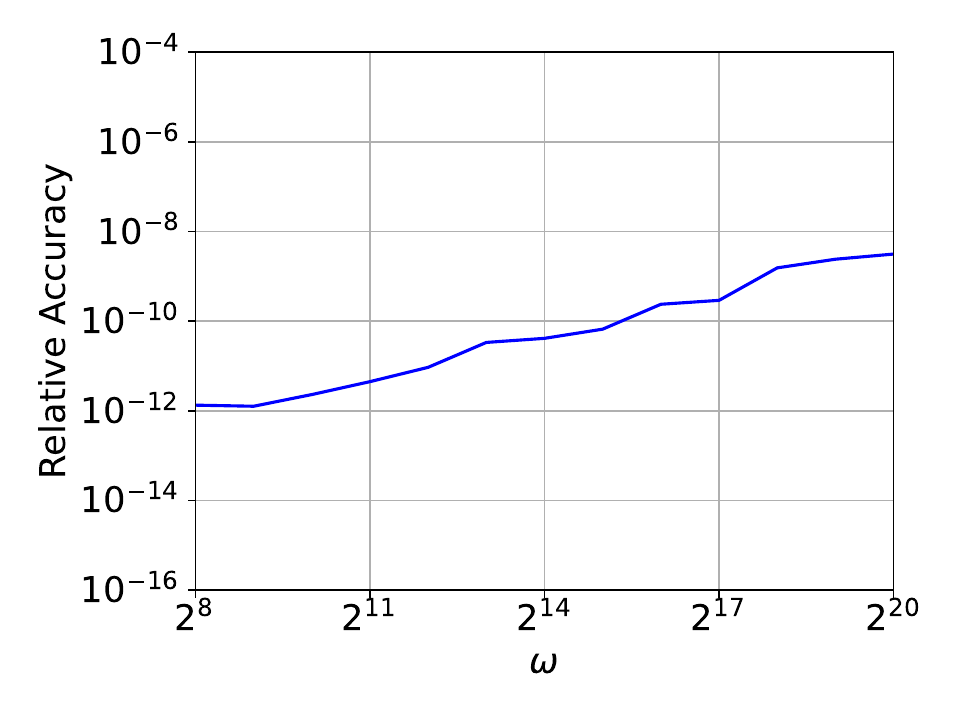}
\hfil

\hfil
\includegraphics[width=.33\textwidth]{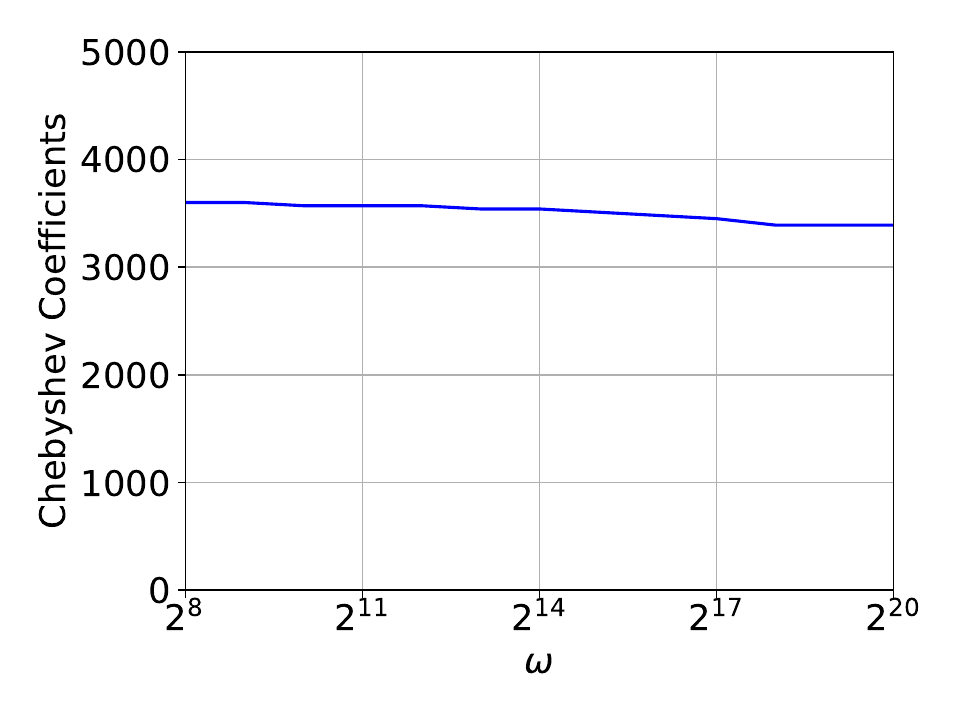}
\hfil
\includegraphics[width=.33\textwidth]{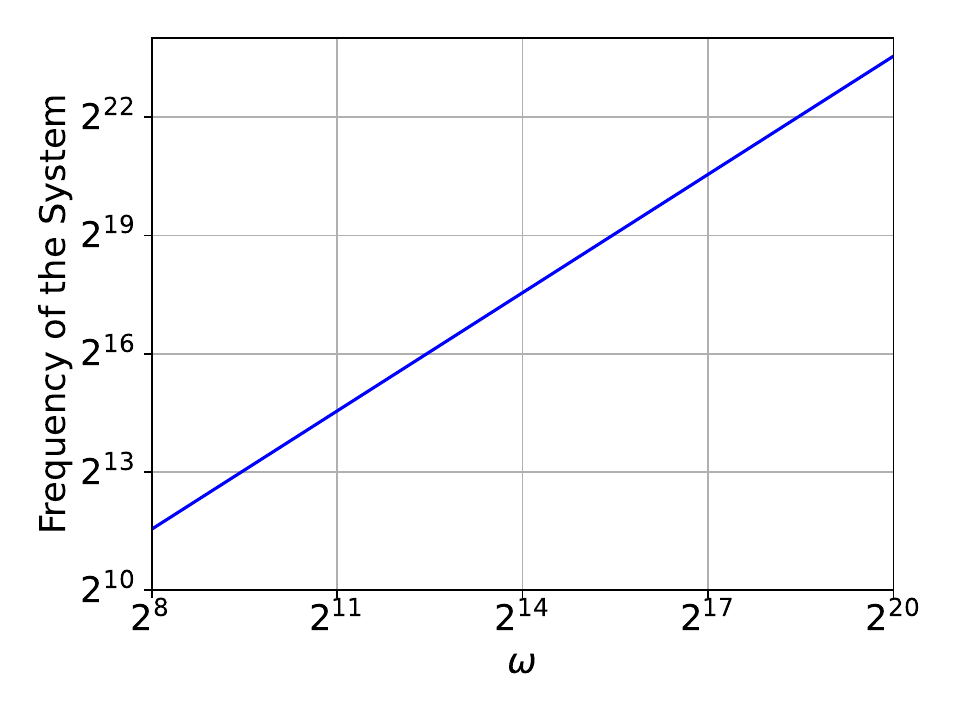}
\hfil

\caption{The results of the experiment of Subsection~\ref{section:experiments:4}.
The plot in the upper left gives the time required by our method
as a function of the parameter $\omega$.  The upper-right plot gives the 
 error in the solution of the initial value problem for the system (\ref{experiments4:system})
as a function of $\omega$.  
The plot on the lower left  shows the total number
of Chebyshev coefficients used to represent the solutions of (\ref{experiments4:system}),
again as a function of the parameter $\omega$.
The plot on the lower right gives the frequency $\Omega$ of the system as a function of
the parameter $\omega$.
}
\label{experiments4:figure1}
\end{figure}

\vfil
\begin{figure}[h!]
\centering

\hfil
\includegraphics[width=.32\textwidth]{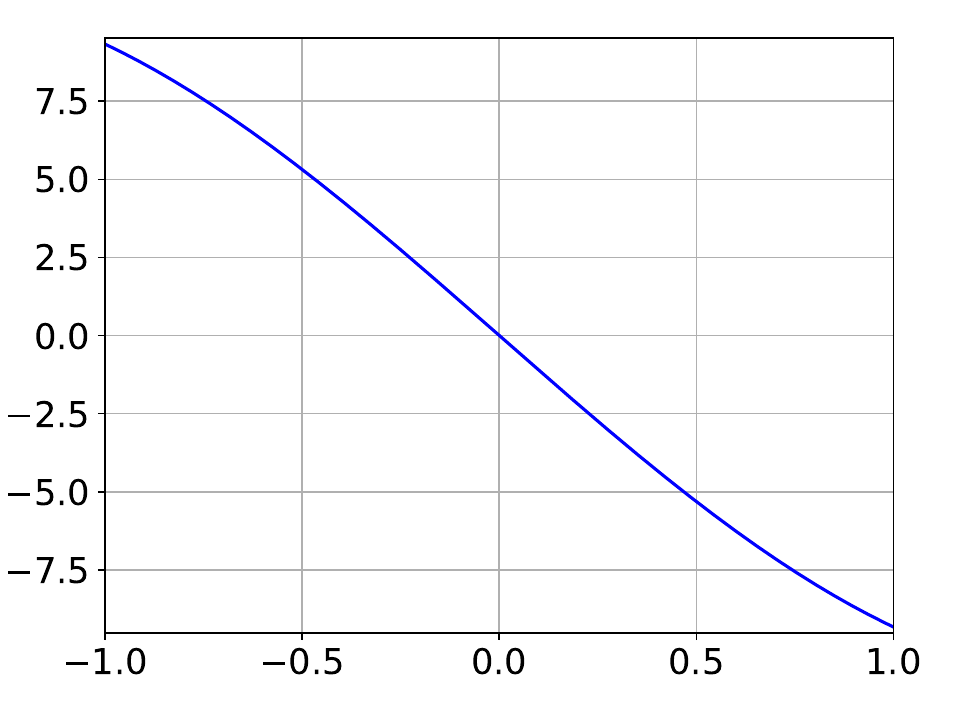}
\hfil
\includegraphics[width=.32\textwidth]{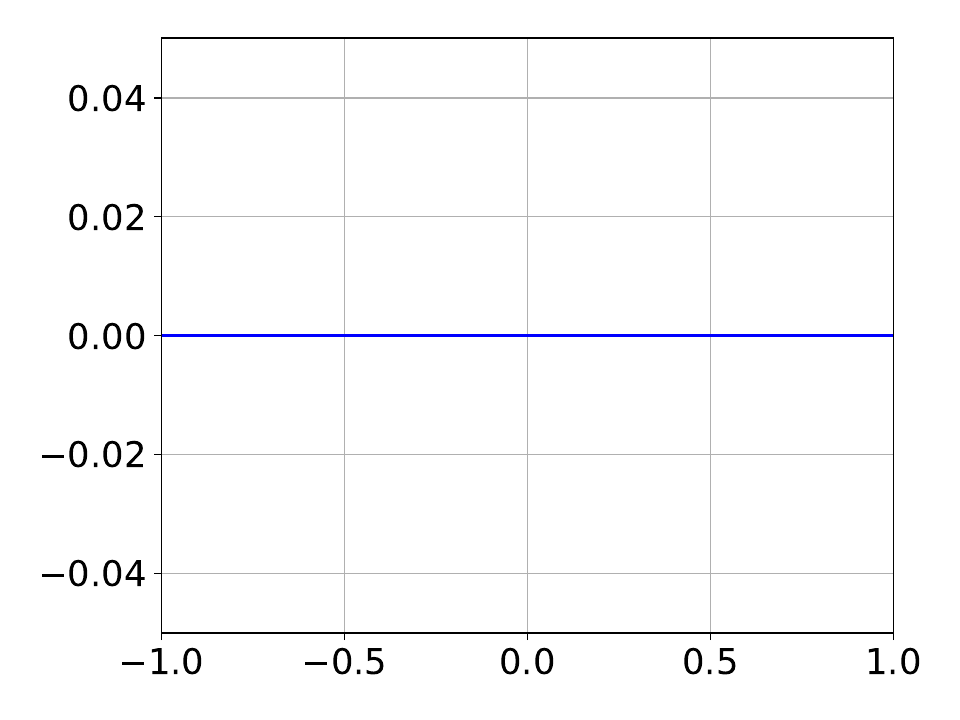}
\hfil
\includegraphics[width=.32\textwidth]{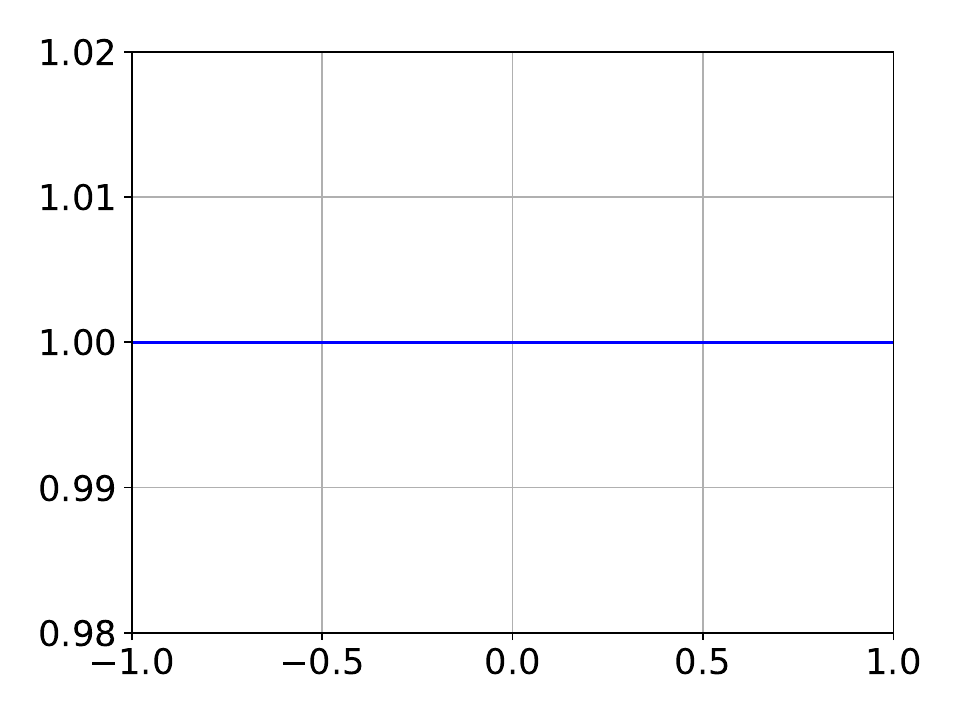}
\hfil

\hfil
\includegraphics[width=.32\textwidth]{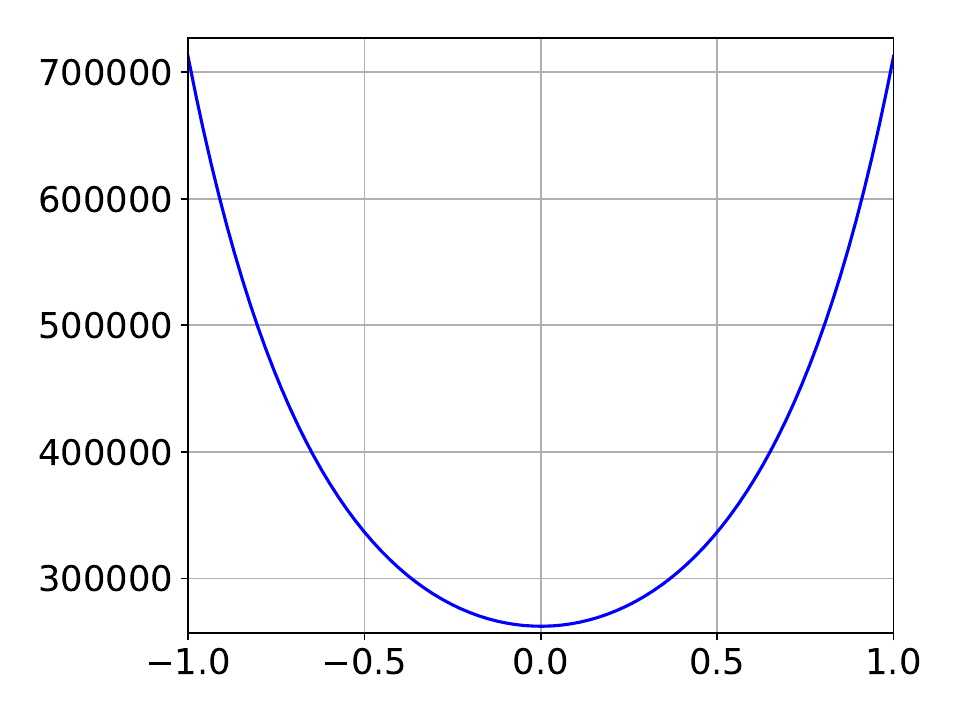}
\hfil
\includegraphics[width=.32\textwidth]{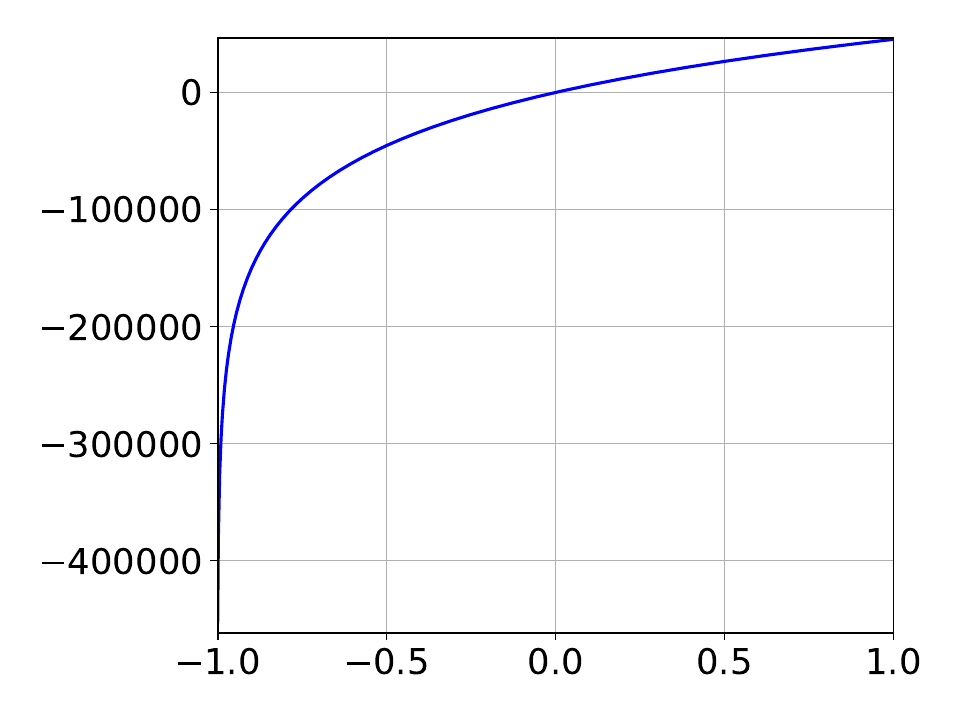}
\hfil
\includegraphics[width=.32\textwidth]{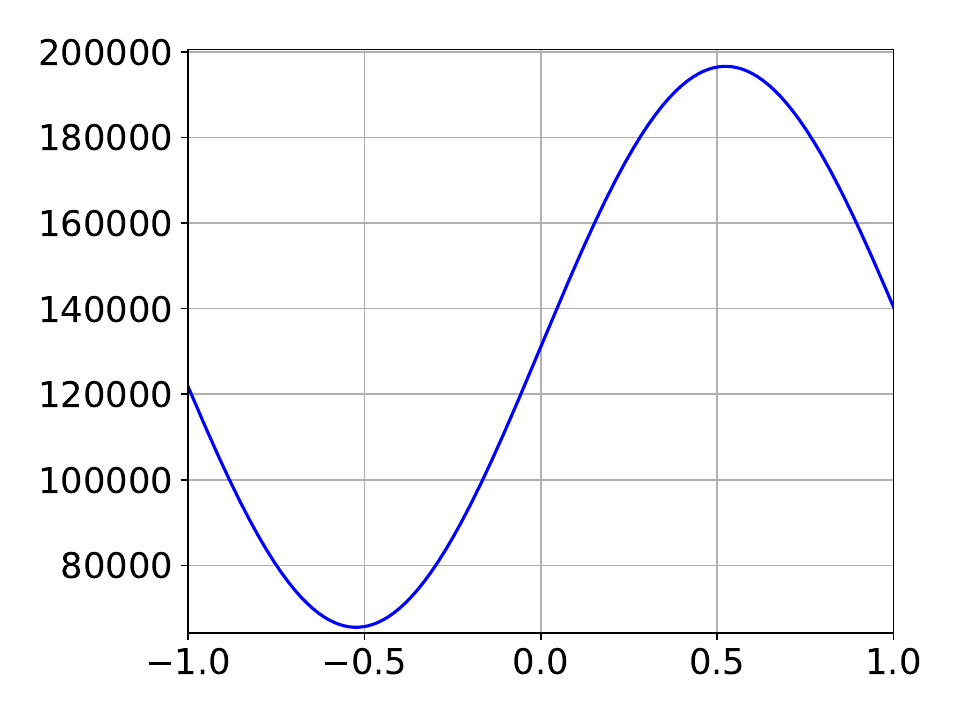}
\hfil

\captionof{figure}{The eigenvalues $\lambda_1(t)$, $\lambda_2(t)$, and $\lambda_3(t)$ of the 
coefficient matrix  $A(t)$
of  Subsection~\ref{section:experiments:4}
when the parameter $\omega$ is equal to $2^{16}$.
Each column corresponds to one eigenvalue,
with the real part appearing in the first row and the imaginary part
in the second. 
}
\label{experiments4:figure2}

\end{figure}

\vfil\eject

Figure~\ref{experiments4:figure1} gives the results of this experiment, while Figure~\ref{experiments4:figure2}
contains plots of the eigenvalues $\lambda_1(t),\lambda_2,\lambda_3(t)$ when $\omega=2^{16}$.
The frequency $\Omega$ of the systems considered ranged from around 2,996 (when $\omega=2^8$) to 
approximately 12,270,000 (when $\omega=2^{20}$).  
The number of Chebyshev coefficients used to represent the solutions of (\ref{experiments4:system})
decreases gradually as $\omega$ increases, from a high of 3,570 ($\omega=2^8$) to a low 
of $3,390$ (when $\omega=2^{18}, 2^{19}, 2^{20}$).  Likewise the running time of the algorithm
trended downward as $\omega$ increases, falling from approximately $29$ milliseconds
when $\omega=2^8$ to roughly $21$ milliseconds when $\omega=2^{20}$.  Again, this is expected
because the complexity of the phase functions typically decreases with increasing frequency.

\end{subsection}

%
%
%
\begin{subsection}{An initial value problem for a system of four equations}
\label{section:experiments:5}

In this experiment, we solved the system of differential equations
\begin{equation}
\mathbf{y}'(t) = A(t) \mathbf{y}(t),
\label{experiments5:system}
\end{equation}
where coefficient matrix $A(t)$ is given by
\begin{equation*}
\left(
\begin{array}{ccc}
 \frac{\cos (t) \left(\log (t+2)-i \sqrt{k}\right)-8 i k e^{2 t^2}}{4 e^{t^2}+\cos (t)} & 0 & -\frac{ie^{t^2} \cos (t) \left(2 k e^{t^2}-\sqrt{k}-i \log (t+2)\right)}{\left(t^2+1\right) \left(4e^{t^2}+\cos (t)\right)} \\
 0 & \frac{-i k (t-8)+t \left(2 e^t-i k\right) \sin (t)+2 e^t t}{4 t^2-2 t-2 t \sin (t)+4} & 0 \\
 -\frac{4 i \left(t^2+1\right) \left(2 k e^{t^2}-\sqrt{k}-i \log (t+2)\right)}{4 e^{t^2}+\cos (t)} & 0& \frac{2 e^{t^2} \left(-i k \cos (t)-2 i \sqrt{k}+2 \log (t+2)\right)}{4 e^{t^2}+\cos (t)} \\
 0 & \frac{t \left(2 e^t \left(t^2+1\right)-i k \left(t^2-3\right)\right)}{2 \left(t^2+1\right) \left(2t^2-t-t \sin (t)+2\right)} & 0 \\
\end{array}
\right.
\end{equation*}

\begin{equation*}
\left.
\begin{array}{c}
 0 \\
 -\frac{\left(2 e^t \left(t^2+1\right)-i k \left(t^2-3\right)\right) (\sin (t)+1)}{2 t^2-t-t \sin(t)+2} \\
 0 \\
 \frac{i k \left(t^4+2 t^2-2 t+1\right)-2 i k t \sin (t)-2 e^t\left(t^2+1\right)^2}{\left(t^2+1\right) \left(2 t^2-t-t \sin (t)+2\right)} \\
\end{array}
\right),
\end{equation*}
over the interval $[-1,1]$ subject to the condition
\begin{equation}
\mathbf{y}(0) = 
\left(
\begin{array}{c}
1\\
-1\\
1\\
-1
\end{array}
\right).
\end{equation}
\vfil\eject

\begin{figure}[h!]
\hfil
\includegraphics[width=.37\textwidth]{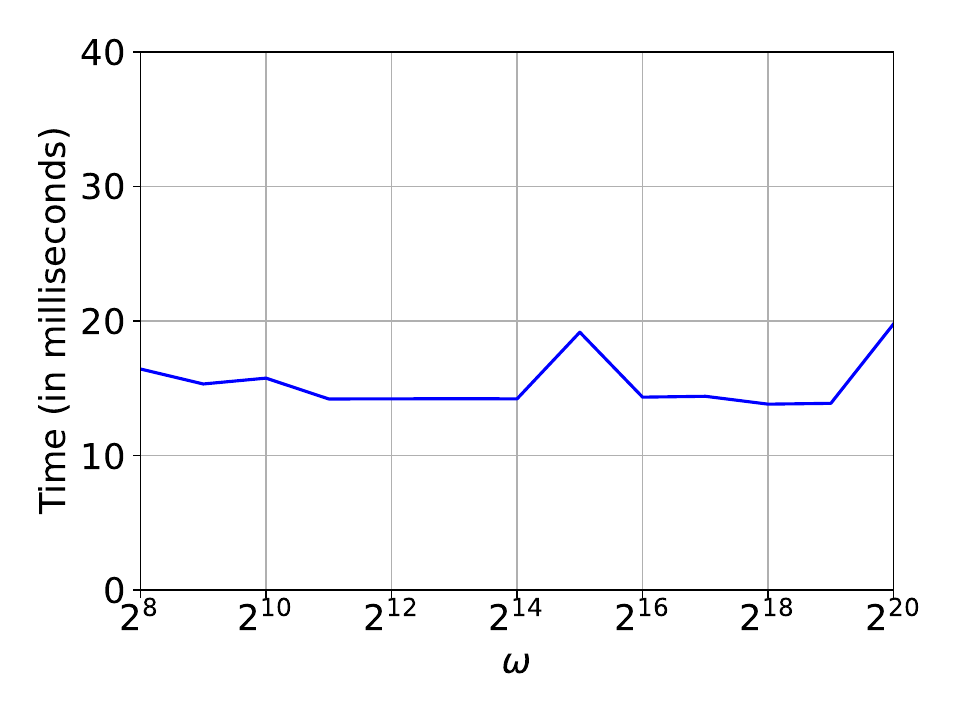}
\hfil
\includegraphics[width=.37\textwidth]{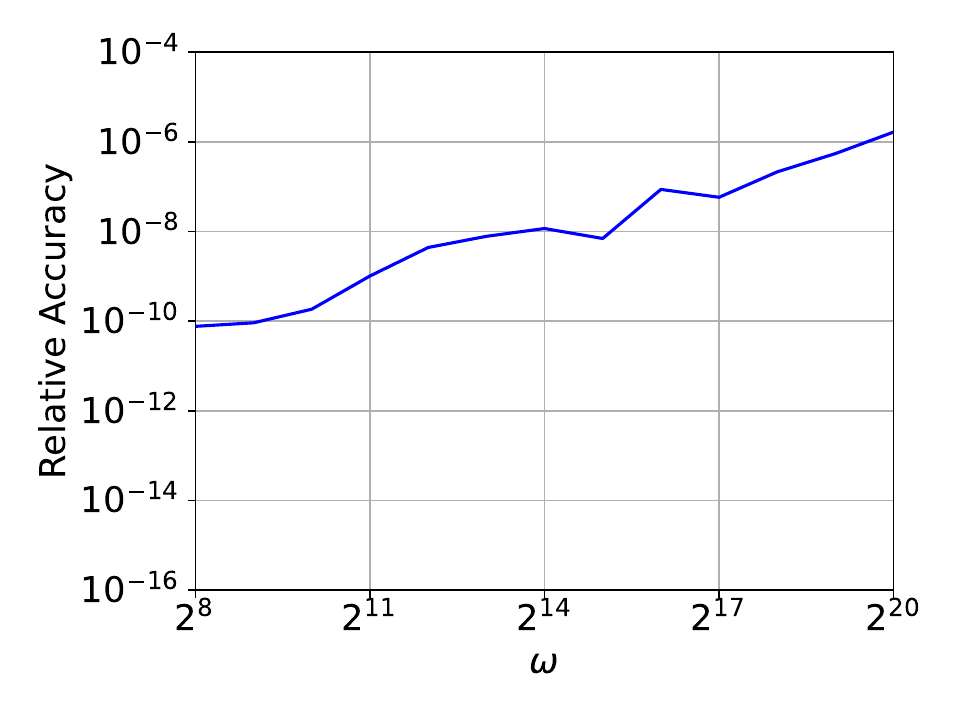}
\hfil

\hfil
\includegraphics[width=.37\textwidth]{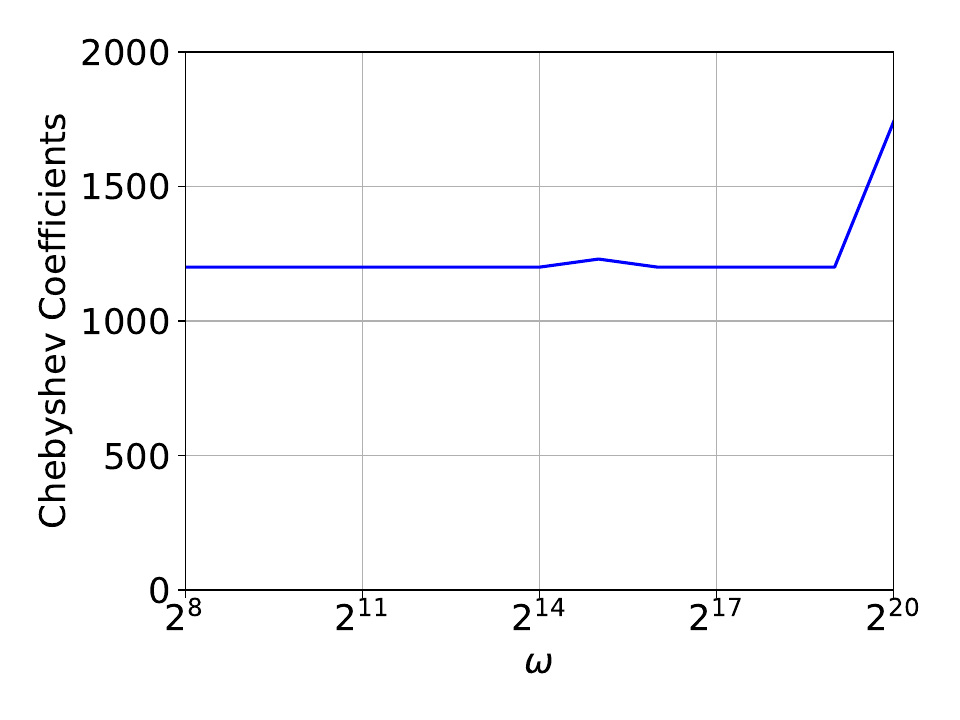}
\hfil
\includegraphics[width=.37\textwidth]{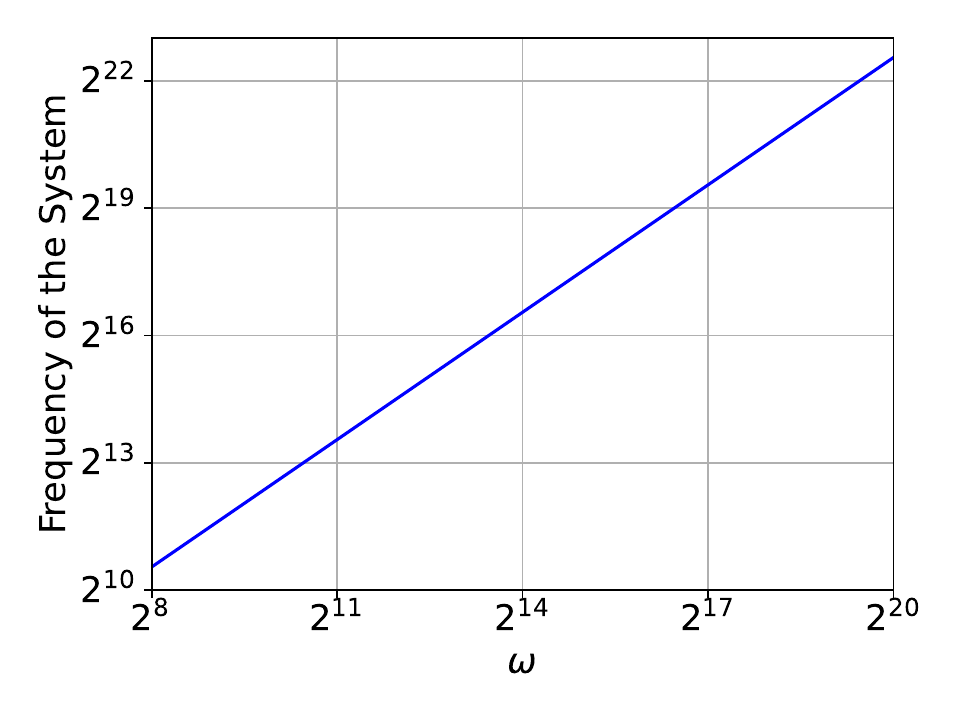}
\hfil

\caption{The results of the experiment of Subsection~\ref{section:experiments:5}.
The plot in the upper left gives the time required by our method
as a function of the parameter $\omega$.  The upper-right plot gives the 
absolute error in the solution of the initial value problem for the system (\ref{experiments5:system})
as a function of $\omega$.  
The plot on the lower left  shows the total number
of Chebyshev coefficients required to represent the solutions of the system (\ref{experiments5:system}), again as a function of the parameter
$\omega$.   The plot on the lower right gives the frequency $\Omega$ of the system as a function of
the parameter $\omega$.}
\label{experiments5:figure1}
\end{figure}

\vfil

\begin{figure}[h!]
\centering

\hfil
\includegraphics[width=.24\textwidth]{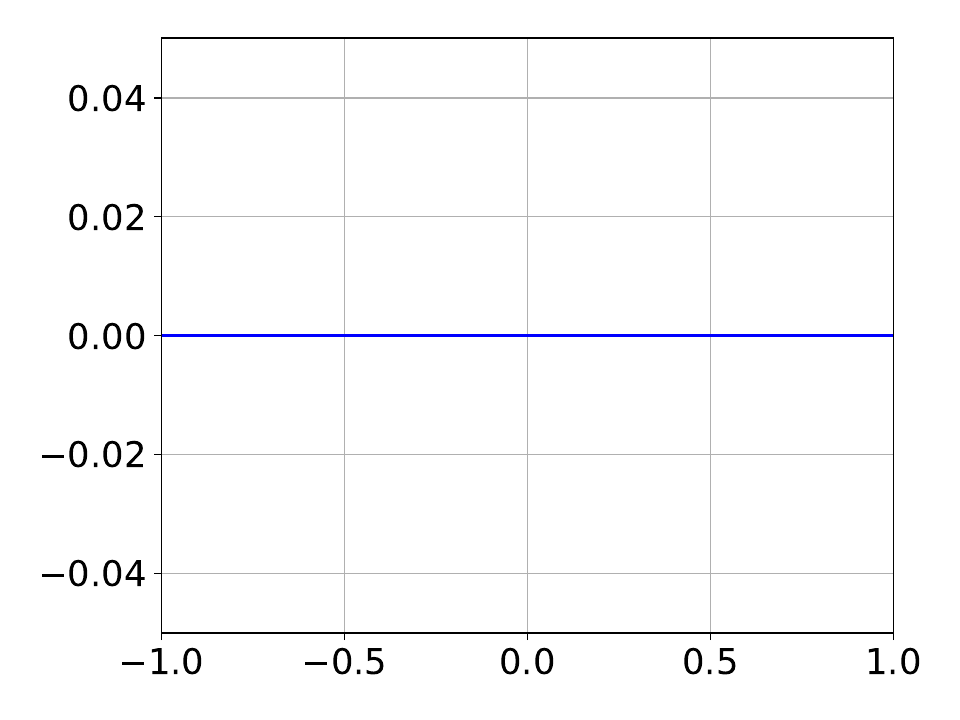}
\hfil
\includegraphics[width=.24\textwidth]{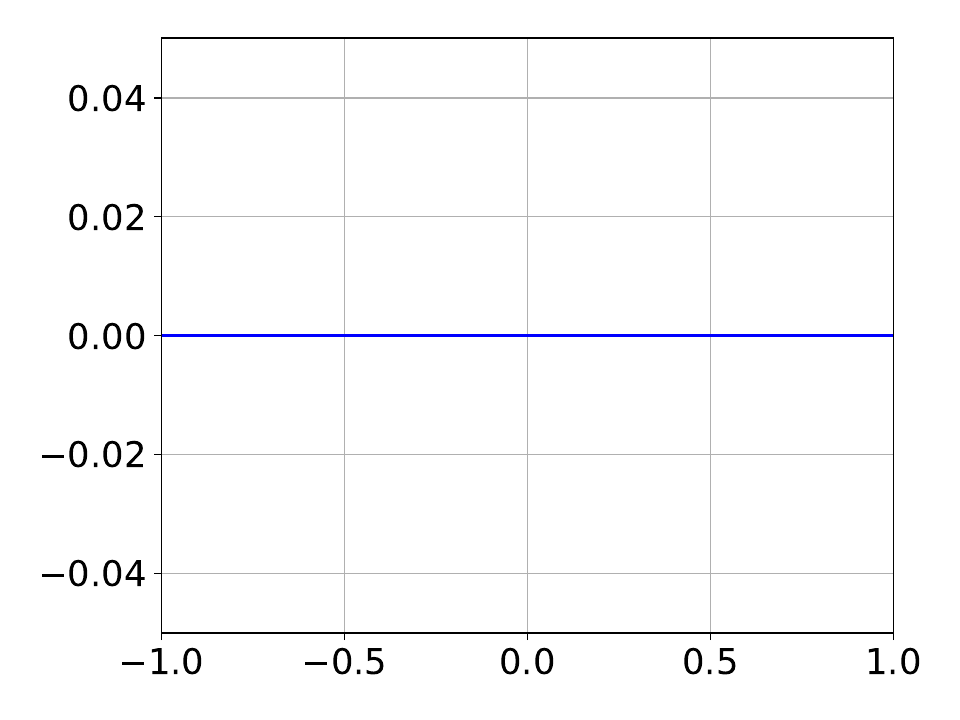}
\hfil
\includegraphics[width=.24\textwidth]{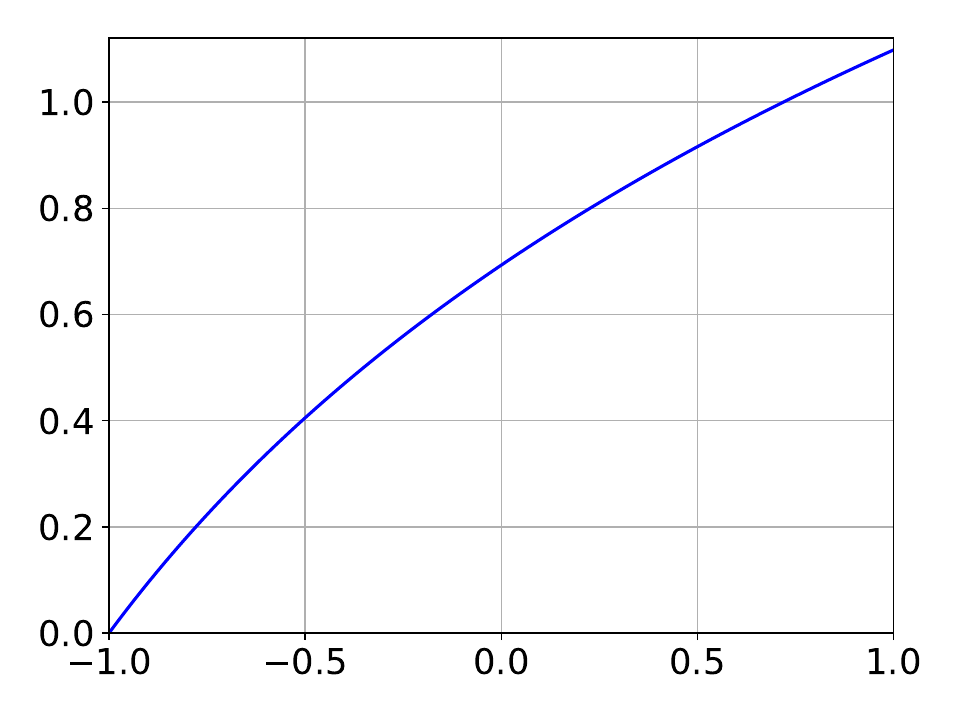}
\hfil
\includegraphics[width=.24\textwidth]{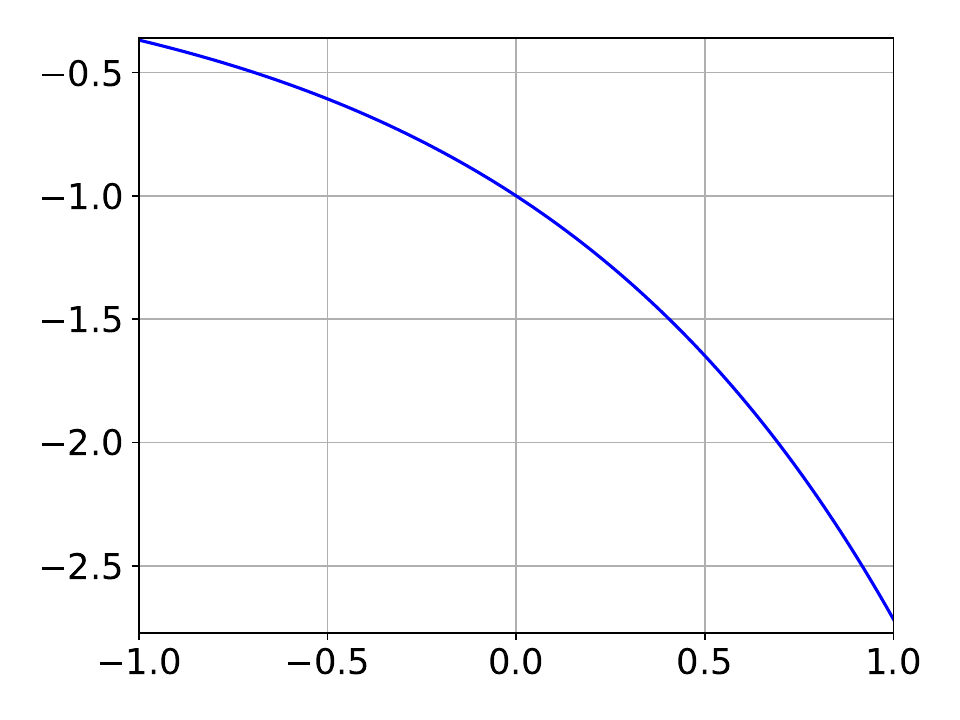}
\hfil

\hfil
\includegraphics[width=.24\textwidth]{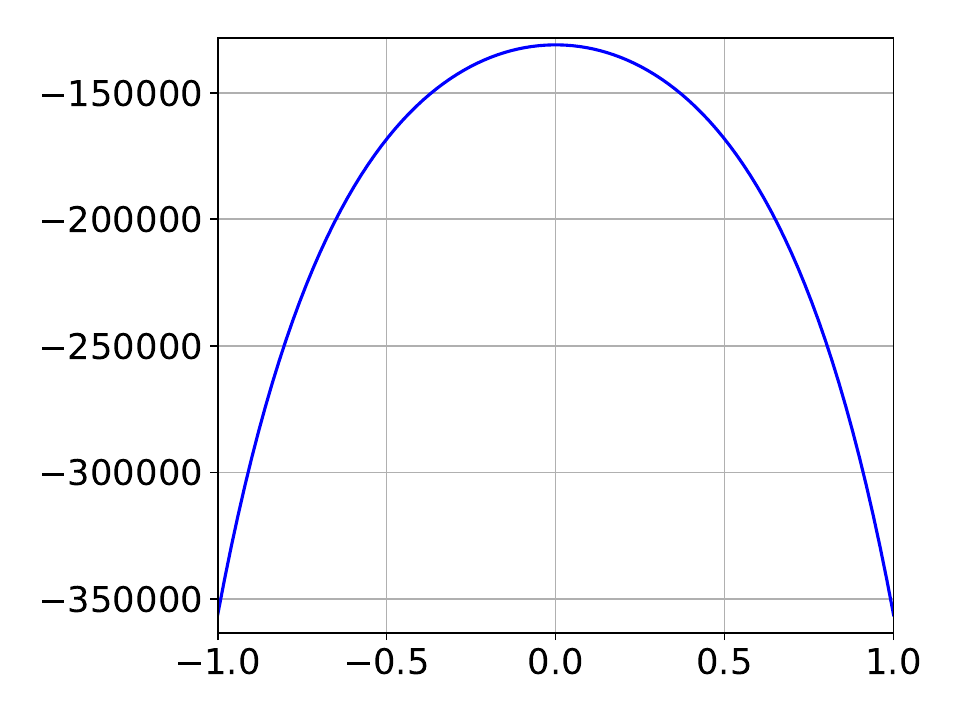}
\hfil
\includegraphics[width=.24\textwidth]{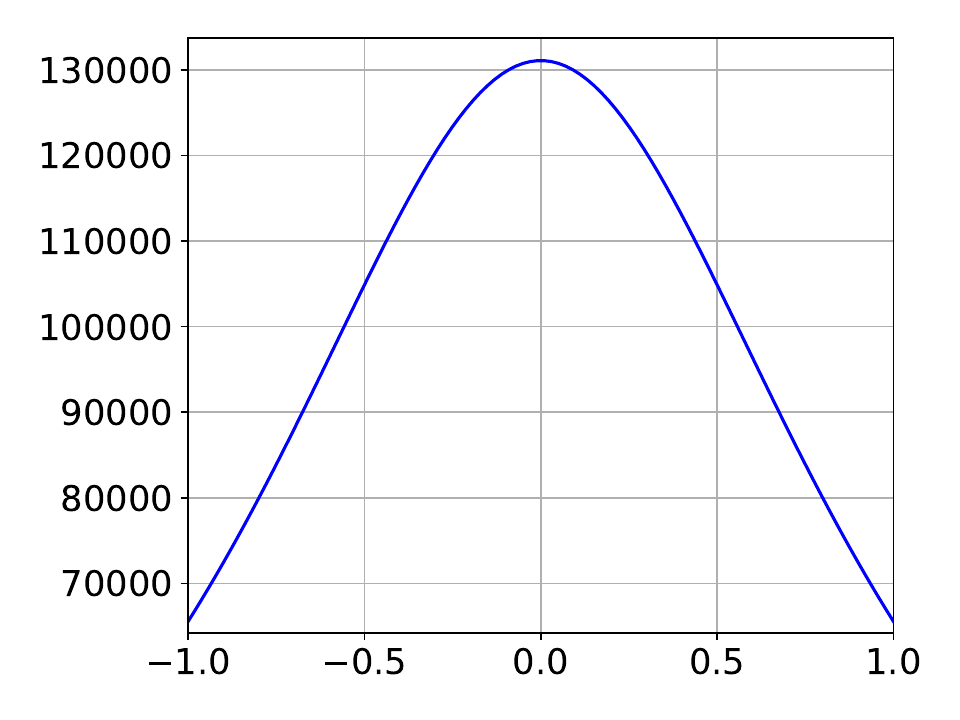}
\hfil
\includegraphics[width=.24\textwidth]{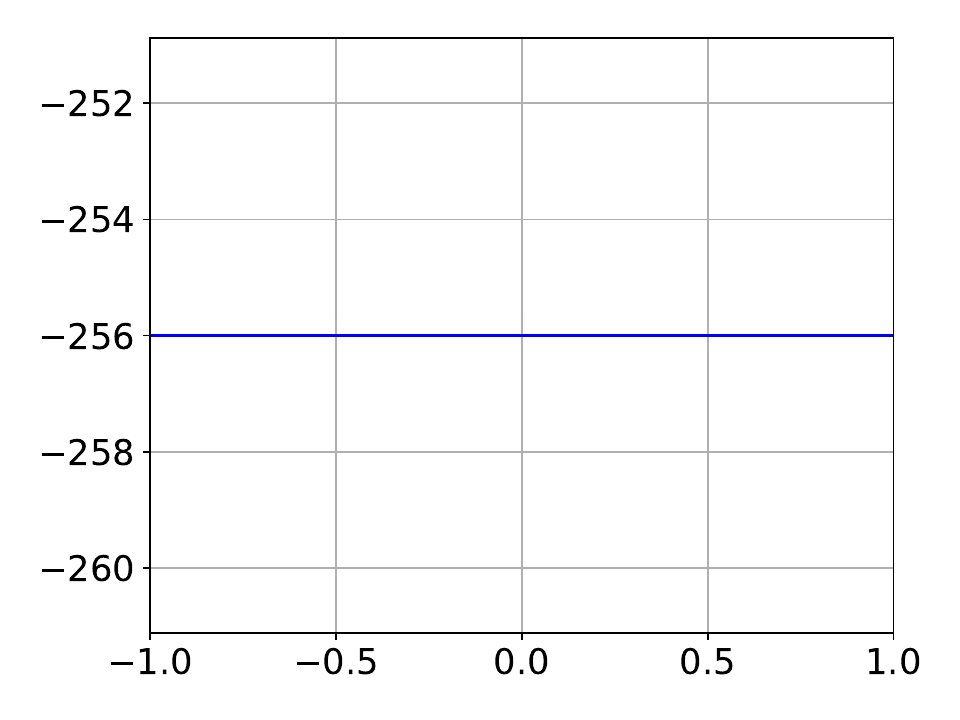}
\hfil
\includegraphics[width=.24\textwidth]{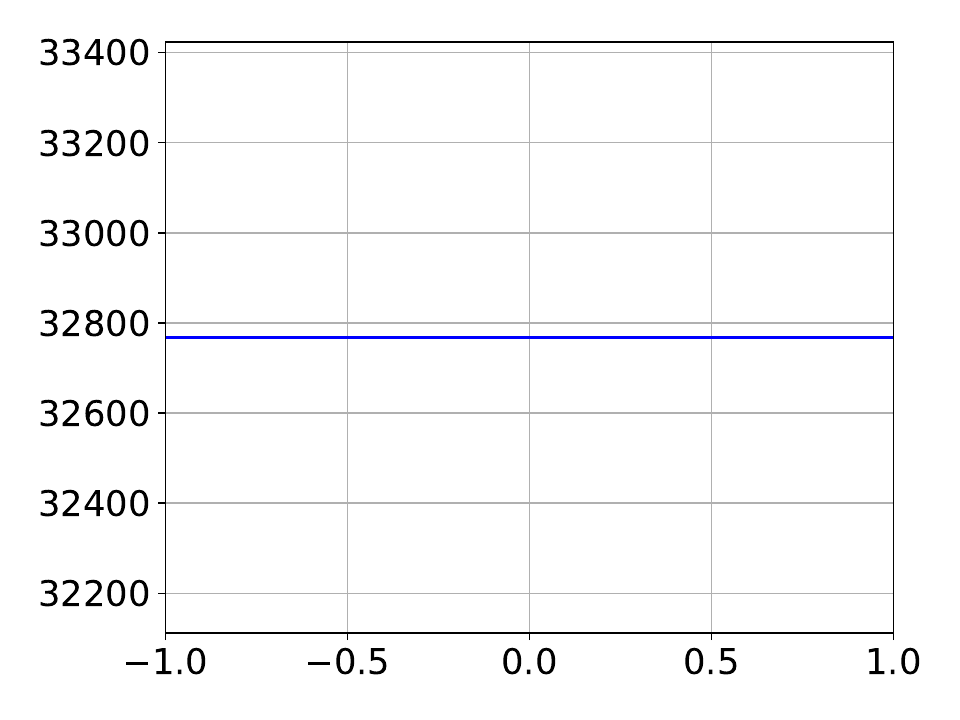}
\hfil

\captionof{figure}{The imaginary parts
eigenvalues $\lambda_1(t), \lambda_2(t), \lambda_3(t), \lambda_4(t)$ of the 
coefficient matrix  $A(t)$
of  Subsection~\ref{section:experiments:5}
when the parameter $\omega$ is equal to $2^{16}$.
Each column corresponds to one eigenvalue,
with the real part appearing in the first row and the imaginary part
in the second. 
}
\label{experiments5:figure2}

\end{figure}

\vfil\eject

The eigenvalues of $A(t)$ are given by the following formulas:
\begin{equation}
\begin{aligned}
\lambda_1(t) &= -2i\omega\exp(t^2),\\
\lambda_2(t) &= \frac{2i\omega }{1+t^2},\\
\lambda_3(t) &= \log(2+t) - i \sqrt{\omega} \ \ \mbox{and} \\
\lambda_4(t) &= -\exp(t) + \frac{i \omega}{2}.
\end{aligned}
\end{equation}

The Levin procedure was performed on the interval $[-0.25,0.00]$, and 
the accuracy parameters  were taken to be 
\begin{equation}
\epsilon_{\mbox{\tiny disc}} = 1.0 \times 10^{-10} \ \ \ \mbox{and}\ \ \ 
\epsilon_{\mbox{\tiny phase}} = 1.0 \times 10^{-10}.
\end{equation}
Because the condition number of the transformations $\Phi$ our algorithm forms tend to increase with the order of the system,
the accuracy achievable by the algorithm of this paper generally decreases with order.  
This is why the accuracy parameters needed to be lowered somewhat for this experiment.
The initial vector for constructing the transformation matrix was taken to be 
\begin{equation}
\mathbf{v} = \left(
\begin{array}{c}
0\\
1\\
1\\
0
\end{array}
\right).
\end{equation}

Figure~\ref{experiments5:figure1} gives the results of this experiment, while Figure~\ref{experiments5:figure2}
contains plots of the eigenvalues $\lambda_1(t)$, $\lambda_2(t)$, $\lambda_3(t)$ and $\lambda_4(t)$ when $\omega=2^{16}$.

\end{subsection}

\end{section}

\begin{section}{Conclusions}
\label{section:conclusions}

We have introduced a numerical method for solving a large class of systems of ordinary differential
equations of modest dimensions in time independent of ``frequency,''  frequency being a measure 
of the magnitudes of the eigenvalues of a system's coefficient matrix.
It operates by transforming the system into a scalar equation via a standard approach, and then
calculating slowly-varying phase functions which represent
a basis in the space of solutions of the scalar equation.  

Our algorithm is limited to systems of modest dimensions because the transformation
matrices we use become increasingly ill-conditioned as the dimension of the system grows.
The set of transformations which take a given system to scalar form, however, is well-known
to be  large in various senses  and it is likely that extremely well-conditioned transformations exist.  The most obvious way to 
overcome the difficulties of our current approach would be to develop an optimization algorithm
to efficiently search through the set of possible transformations for a well-conditioned  one.     
The development of a numerical version of the elimination algorithm of \cite{Barkatou} 
which we briefly discuss in Section~\ref{section:reduction}
is another possible approach to constructing a stable transformation matrix.
Alternatively, one could simply insert the representation 
(\ref{introduction:expdiag}) into the system (\ref{introduction:system}) and derive a set of differential equations 
for the phase functions $\psi_1,\ldots,\psi_n$ and transformation matrix $\Phi(t)$.
The slowly-varying phase functions are effectively uniquely determined as long as the frequency
of the system is large enough (see \cite{BremerRokhlin}),  but the choice of the  transformation 
$\Phi(t)$ is  highly nonunique, which complicates such an approach.
The authors are actively investigating these and other possible methods for improving the 
algorithm of this paper.

Despire its limitations, the method described here clearly has many applications in scientific problems.
Moreover, we view it a promising  first step toward the development of robust, high-accuracy frequency-independent
solvers for systems of ordinary differential equations.

\end{section}

\begin{section}{Acknowledgments}
JB  was supported in part by NSERC Discovery grant  RGPIN-2021-02613.
The authors are grateful to Kirill Serkh for many helpful discussions regarding
this work.
\end{section}

\begin{section}{Data availability statement}
The datasets generated during and/or analysed during the current study are available from the corresponding author on reasonable request.
\end{section}

\bibliographystyle{acm}
\bibliography{scalar.bib}


\end{document}